\documentclass[11pt]{article}
\usepackage[utf8]{inputenc}
\usepackage{algpseudocode}
\usepackage{xcolor}
\usepackage{amsmath}
\usepackage{amssymb}
\usepackage{siunitx}
\usepackage{float}
\usepackage{mathtools}
\usepackage[numbers,sort&compress]{natbib}
\usepackage[edges]{forest}
\usepackage{setspace}
\usepackage{algorithm2e}
\SetKwComment{Comment}{/* }{ */}
\RestyleAlgo{ruled}
\DontPrintSemicolon

\usepackage{doi}

\usepackage{authblk}
\usepackage[symbol]{footmisc}

\usepackage{orcidlink}

\usepackage[english]{babel}
\usepackage[font=small,labelfont=bf]{caption}

\usepackage[skip=8mm]{caption} 
\usepackage{subcaption}
\usepackage{graphicx}

\usepackage[short,nocomma,c1]{optidef}

\usepackage{duckuments}

\usepackage{tikz}
\usetikzlibrary{arrows,calc,fit,patterns}
\usepackage{eurosym}

\usepackage{pdfcomment}

\usepackage{pgfplots,pgfplotstable}
\usepgfplotslibrary{fillbetween}
\usepackage{enumerate}
\usepackage[shortlabels]{enumitem}

\usepackage{hyperref}
\usepackage{amsthm}
\usepackage[capitalize,nameinlink]{cleveref}

\usepackage{tabularx}

\usepackage[title]{appendix}

\SetKwComment{Comment}{\# }{} 

\newcommand{\crefdefpart}[2]{%
	\hyperref[#2]{\namecref{#1}~\labelcref*{#1}~\ref*{#2}}%
}

\usepackage[left=2cm,top=1.4cm,right=2cm,bottom=2.5cm,nohead]{geometry}
\usepackage{amssymb}
\usepackage{amsmath}
\usepackage{bbold}

\makeatletter
\renewcommand{\section}{\@startsection{section}{1}{0ex}%
                                   {-3.5ex \@plus -1ex \@minus -.2ex}%
                                   {0.1ex \@plus.2ex}%
                                   {\normalfont\Large\bfseries\sffamily}}
\renewcommand{\subsection}{\@startsection{subsection}{2}{0ex}%
                                     {-3.25ex\@plus -1ex \@minus -.2ex}%
                                     {1.5ex \@plus .2ex}%
                                     {\normalfont\large\bfseries\sffamily}}
\renewcommand{\subsubsection}{\@startsection{subsubsection}{3}{0ex}%
                                     {-3.25ex\@plus -1ex \@minus -.2ex}%
                                     {1.5ex \@plus .2ex}%
                                     {\normalfont\normalsize\bfseries\sffamily}}
\renewcommand{\paragraph}{\@startsection{paragraph}{4}{\z@}%
                                    {3.25ex \@plus1ex \@minus.2ex}%
                                    {-1em}%
                                    {\normalfont\normalsize\bfseries\sffamily}}
\renewcommand{\subparagraph}{\@startsection{subparagraph}{5}{\parindent}%
                                       {3.25ex \@plus1ex \@minus .2ex}%
                                       {-1em}%
                                      {\normalfont\normalsize\bfseries\sffamily}}
\renewcommand{\@maketitle}{%
  \newpage
  \null
  \begin{center}%
  \let \footnote \thanks
    {\Large \textsf{\textbf{\@title}} \par}%
    \vskip 0.5em%
    {\large
      \lineskip .5em%
      \begin{tabular}[t]{c}%
        \textsl{\@author}
      \end{tabular}\par}
  \end{center}%
  \par
  \vskip 1em}
\makeatother

\usepackage{tikz}

\title{Towards net-zero manufacturing: carbon-aware scheduling for GHG emissions reduction}

\author[,a,b]{Andrea Mencaroni\,\orcidlink{0000-0002-0110-3218}\thanks{Corresponding author. \textit{E-mail}: \href{mailto:andrea.mencaroni@ugent.be}{andrea.mencaroni@ugent.be}}}
\author[a,b]{Pieter Leyman}
\author[a]{Birger Raa}
\author[a,b]{Stijn De Vuyst}
\author[a,b]{Dieter Claeys}
\affil[a]{\footnotesize Department of Industrial Systems Engineering and Product Design, Ghent University, Ghent,
Belgium} 
\affil[b]{\footnotesize Industrial Systems Engineering (ISyE), Flanders Make, Kortrijk, Belgium}
\date{} 

\providecommand{\keywords}[1]{\textbf{\textit{Keywords --}} #1}

\usepackage{xpatch}
\makeatletter
\AtBeginDocument{\xpatchcmd{\@thm}{\thm@headpunct{.}}{\thm@headpunct{}}{}{}}
\makeatother

\begin{document}
\maketitle

\section*{Abstract}
Detailed scheduling has traditionally been optimized for the reduction of makespan and manufacturing costs.
However, growing awareness of environmental concerns and increasingly stringent regulations are pushing manufacturing towards reducing the carbon footprint of its operations.
Scope 2 emissions, which are the indirect emissions related to the production and consumption of grid electricity, are in fact estimated to be responsible for more than one-third of the global GHG emissions.
In this context, carbon-aware scheduling offers a promising strategy for lowering manufacturing’s carbon footprint by accounting for the time-dependent carbon intensity of the grid and the availability of on-site renewable electricity.

This study introduces a carbon-aware permutation flow-shop scheduling model aimed at reducing scope 2 emissions.
The model is formulated as a mixed-integer linear problem, incorporating forecasted grid generation mix, available on-site renewable electricity, and job-specific power requirements.
To solve the problem, we propose a dual random-key memetic algorithm that integrates evolutionary strategy with local search.

Computational experiments show that substantial reductions in carbon emissions can be achieved through carbon-aware scheduling, with limited impact on makespan.
For instance, emission reductions of up to $47.6 \, \%$ are observed with only a $1.8 \, \%$ increase in makespan in single-machine scenarios.
These results demonstrate the potential of carbon-aware scheduling to leverage time-dependent energy-related data and reduce GHG emissions without severely compromising operational performance.

\keywords{Scheduling, Carbon-aware, Sustainability, Evolutionary computing, Permutation flow-shop}

\section{Introduction}
\label{section:introduction}
Amid the global push towards net-zero greenhouse gas (GHG) emissions \cite{eu_greendeal2023}, the energy supply sector, being the largest contributor with over one third of total GHGs \cite{IPCC_AR6}, is set for an unprecedented transformation.
While the ever-growing deployment of renewable energy sources is a crucial step in the EU's path to net-zero \cite{net-zero_plan_2019}, it alone is not sufficient.
The inherent variability in renewable energy availability necessitates addressing the demand side as well, making a shift towards supply-driven energy consumption inevitable \cite{EEA_report_2023}.

Energy is fundamental to all sectors of the economy, from transport to households, agriculture, and industry.
With the latter accounting for roughly one fourth of total EU energy consumption in 2022 \cite{Eurostat_energy_report_2023}, how energy is generated becomes particularly critical for emissions reduction.
GHG emissions are categorized in three scopes: scope 1, which regards direct emissions from owned or controlled sources; scope 2, covering indirect emissions from the generation of purchased electricity, and scope 3, which encompasses all other indirect emissions from an organization's value chain (i.e., its upstream and downstream activities) \cite{GHGCorporateStandard}.
The mix of energy sources used to produce purchased electricity directly affect scope 2 emissions and, consequently, has a major impact on a company's overall carbon footprint.

In a typical manufacturing environment, electricity is supplied by a combination of on-site renewable utility systems and the public grid \cite{BAUMGARTNER2019755}.
While using on-site renewables is preferable due to their zero scope 2 emissions \cite{GHGscope2guidance}, their availability is limited and subject to variability \cite{Batalla-Bejerano2016}.
In contrast, grid electricity can be considered effectively unlimited in capacity because external balancing mechanisms ensure that supply matches demand. This is achieved by adjusting the generation mix -- i.e., varying the contribution of different energy sources such as natural gas, coal, and nuclear power -- to meet the real-time electricity needs \cite{Ahlqvist2022}.
As a result, grid electricity's scope 2 emissions vary over time, leading to time-dependent carbon intensity \cite{Miller2022}.
Therefore, there is substantial potential to reduce GHG emissions by aligning electricity consumption with the availability of on-site renewable electricity and the grid's carbon intensity \cite{KOPSAKANGASSAVOLAINEN2017384}.

Hence, the following question naturally arises: \textit{how} can industry optimally respond to fluctuations in grid carbon intensity and renewable electricity availability, thereby decreasing GHG emissions?
Carbon-aware scheduling emerges as a promising answer to this research question.

\subsection{Overview of the state-of-the-art}
\label{subsection:overview_state_of_the_art}
Scheduling refers to the well-known decision-making process of allocating resources to tasks over given time periods to optimize one or more objectives \cite{PinedoScheduling}.
Given the complexity of the process, scheduling has traditionally been modeled as an optimization problem subject to constraints.
Many different formulations of scheduling problems exist, reflecting the resource configuration and the nature of the tasks to be scheduled \cite{Baker2018}.
One of the most common formulations is the permutation flow-shop scheduling problem (PFSP) \cite{Fernandez-Viagas2017}.
With a PFSP model, machines are arranged in series, and the sequence of jobs on machines is maintained across all the machines.
The popularity of PFSP lies in the fact that this configuration resembles how assembly lines are built, where semi-finished products are moved through subsequent workstations by means of conveyor systems.
Since the 1960s, the practical relevance of PFSP has drawn considerable research attention.
However, the predominant focus has been on minimizing makespan \cite{Ostermeier2024, Ahmadian2021}.
While optimizing for faster production has driven industrial development for decades, this objective no longer fully aligns with the current context of environmental concerns.
The global energy crisis of 2021, which triggered record-high energy prices \cite{IEA2022}, accelerated a shift in research focus toward minimizing energy costs by incorporating dynamic electricity pricing.
While many studies have addressed the reduction of energy costs, often in combination with makespan minimization, the majority have relied on fixed time-of-use (TOU) pricing strategies (e.g., \cite{Che2016, Cheng2017, Cui2021, Ghorbanzadeh2023, Ho2021, Ho2022, Jiang2020, Karimi2021, Luo2013, Moon2014, Moon2013, Oukil2022, Rubaiee2019, Shen2021, Shen2023, Trevinomartinez2022Neutralityoptimization, Wang2018}).
Few studies considered real-time pricing (RTP), where electricity prices vary throughout the day without predefined tariffs (e.g., \cite{Ding2016Parallel, FazliKhalaf2018, Tian2024}).
While optimization under dynamic energy prices often leads to reduced peaks in demand, this does not necessarily translate to lower emissions \cite{Holland2008}.
Although electricity prices are generally correlated with grid carbon intensity, optimizing solely for cost fails to fully leverage the emissions reduction potential of demand-side flexibility \cite{Gabrek2025}.

Some studies have considered simultaneous minimization of makespan and carbon emissions.
However, the majority relied on annual average carbon emission factors (e.g., \cite{Ding2016Carbon, Dong2022, Fallahi2023, GhorbaniSaber2022, Liu2016, Liu2017}), neglecting the time-dependent variations in grid carbon intensity.
Only very few have addressed the time-dependent nature of grid carbon intensity.
For instance, Zhang et al.\ \cite{Zhang2014} proposed a flow-shop scheduling model to minimize total energy costs and carbon emissions by incorporating dynamic grid carbon intensity and TOU energy prices.
Similarly, Kelley et al.\ \cite{Kelley2018} developed a single-machine scheduling model to minimize GHG emissions using hourly grid generation mix data.
Lastly, Trevino-Martinez et al.\ \cite{Trevinomartinez2022Footprintoptimization} introduced a single-machine scheduling model to jointly minimize total energy and carbon emissions costs by incorporating a carbon tax.
However, all these studies omitted on-site renewable generation and did not present an algorithm capable of solving real-world sized instances.

More broadly, scheduling problems with time-dependent objective functions have been studied in contexts such as time-dependent processing times \cite{Alidaee1999, Cheng2004}, setup times \cite{Stecco2008}, and discrepancy times \cite{Jaehn2016}.
These studies typically involve a single time-dependent factor, often following a predictable pattern.
In contrast, carbon-aware scheduling introduces three distinct time-dependent elements: the grid's carbon intensity, the availability of on-site renewable electricity, and the power requirements of individual jobs.
On-site renewable generation, additionally, is also limited in supply.
This added complexity presents new modeling and optimization challenges that existing scheduling formulations do not address.

To the best of our knowledge, no existing work has proposed a scheduling framework that simultaneously accounts for time-dependent grid carbon intensity, fluctuating on-site renewable generation, and job-specific power requirements.
As a result, there is a clear need for practical, scalable carbon-aware scheduling algorithms capable of effectively minimizing GHG emissions in manufacturing.

\subsection{Scope of the paper}
\label{subsection:scope_of_the_paper}
To address the identified research gap, this study proposes a permutation flow-shop carbon-aware scheduling model aimed at minimizing scope 2 GHG emissions in manufacturing.
The problem is formulated as a Mixed-Integer Linear Program (MILP) and solved using a dual random-key memetic algorithm, which integrates evolutionary computing with local search.

The key contributions of this study are three-fold.
First, it develops a carbon-aware scheduling model for the PFSP, a common abstraction of sequential assembly lines.
The model integrates time-dependent grid carbon intensity, on-site renewable electricity generation, and job-specific power requirements.
Second, it proposes a dedicated memetic algorithm that efficiently finds high-quality solutions for real-world-sized instances within short computation times, making it suitable for practical applications.
Third, it demonstrates through computational experiments that substantial carbon emission reductions do not necessarily come at the cost of increased makespan, especially in instances with lower scheduling flexibility and limited slack time.

The remainder of this paper is structured as follows.
In Section \ref{section:model_formulation}, the model and the associated terminology are introduced.
Then, in Section \ref{section:memetic_algorithm_framework}, the components and parameters of the memetic algorithm are described.
Next, in Section \ref{section:computational_experiments}, computational experiments are presented and discussed.
Finally, conclusions and directions for future research are outlined in Section \ref{section:conclusion}.

\section{Model formulation}
\label{section:model_formulation}

We now introduce the carbon-aware scheduling model by following two steps.
In Section \ref{subsection:electricity_demand_and_grid_carbon_intensity},we first introduce the concept of time-dependent job power requirements and present the methodology for computing grid carbon intensity.
Then, in Section \ref{subsection:CAS_model}, we present the mathematical formulation of the MILP model, including the relevant terminology and underlying assumptions.

\subsection{Job power requirements and grid carbon intensity}
\label{subsection:electricity_demand_and_grid_carbon_intensity}
\begin{figure}[H]
	\centering
	\includegraphics[width=.6\linewidth]{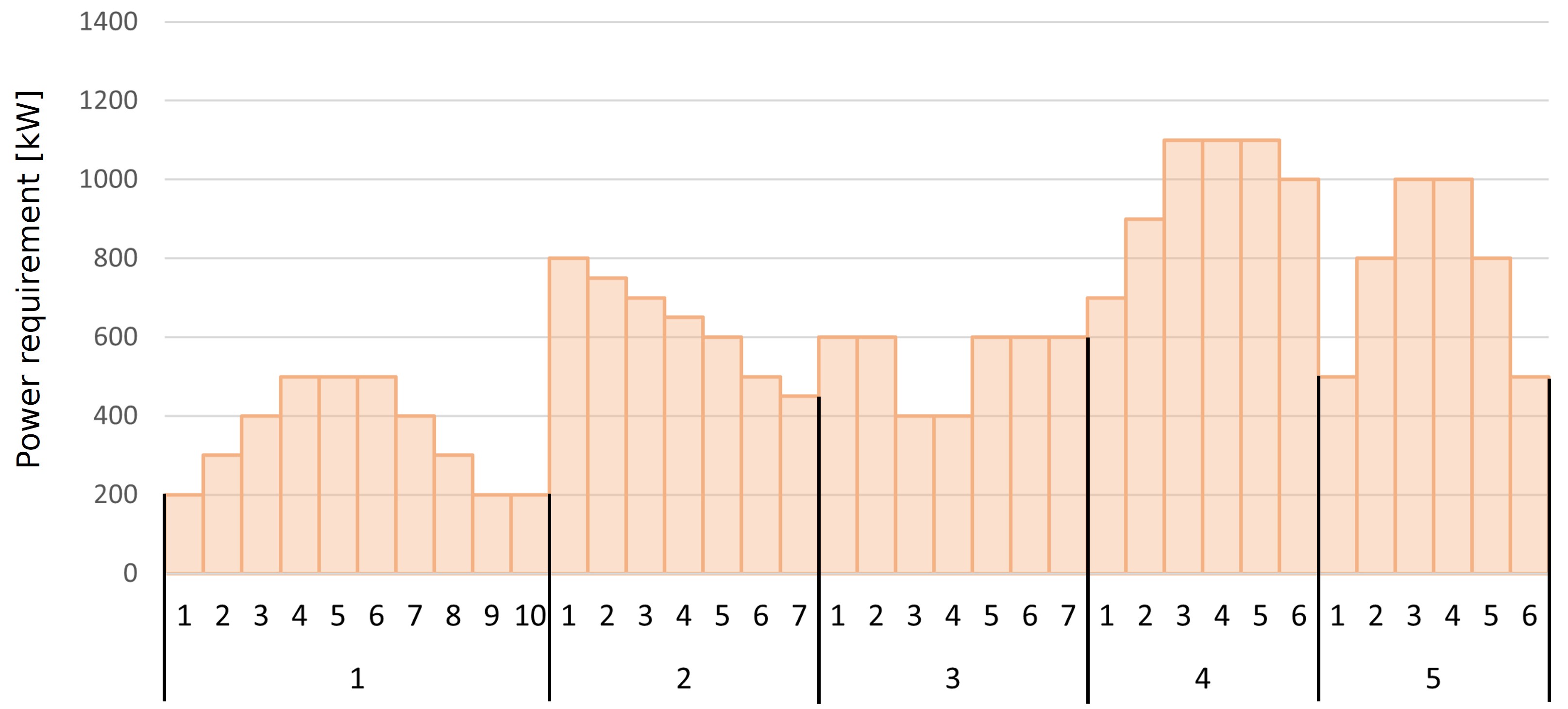}
	\caption{Example of five jobs with corresponding processing time and power requirements}
	\label{fig:job_power_requirements}
\end{figure}

Traditionally, time has been considered the primary resource required to process jobs on machines.
However, manufacturing operations also depend on electrical power, whose demand may vary throughout a job's execution.
For instance, certain tasks may involve a machine warm-up phase, typically resulting in ramp-up or ramp-down power profiles.
Figure \ref{fig:job_power_requirements} illustrates five example jobs with their corresponding processing times and power requirement profiles.

The source of electricity used to power machines significantly impacts GHG emissions.
In this model, two possible sources of electrical power are considered: on-site renewable generation and the public grid.

On-site renewable power includes installations such as photovoltaic (PV) systems, wind turbines, and heat recovery stations.
Electricity generated from these sources is preferable as it does not contribute to scope 2 GHG emissions.
However, on-site renewable power is typically limited in capacity and fluctuates throughout the day depending on external factors such as the weather conditions.

The second source is electricity from the public grid.
Unlike on-site renewable power, the grid provides a virtually unlimited supply, supported by external balancing systems that ensure demand is met at any moment of the day.
However, electricity from the grid is associated with scope 2 GHG emissions, as its production often relies on carbon-intensive energy sources.
These emissions are typically expressed in grams of carbon dioxide equivalent (gCO\textsubscript{2}eq), a standardized unit that accounts for the global-warming potential (GWP) of different greenhouse gasses.

The GHG emissions per unit of grid electricity, which are referred to as grid carbon intensity, are estimated based on the generation mix and the lifecycle emissions of each contributing energy source, which are summarized in Table \ref{tab:lifecycle_emissions_per_source}.
\begin{table}[H]
	\fontsize{9pt}{14pt}\selectfont
	\centering
    \caption{Emissions of selected electricity supply technologies, reprinted from \cite{IPCC_ar5}}
	\begin{tabular}{l|c}
		Source & Lifecycle emissions in gCO\textsubscript{2}eq/kWh (Min/Median/Max)\\
		\hline
		Pulverized coal & 740/\textbf{820}/910 \\
		Gas -- Combined cycle & 410/\textbf{490}/650 \\
		Biomass -- cofiring & 620/\textbf{740}/890 \\
		Biomass -- dedicated & 130/\textbf{230}/420 \\
		Geothermal & 6/\textbf{38}/79 \\
		Hydropower & 1/\textbf{24}/2200 \\
		Nuclear & 3.7/\textbf{12}/110 \\
		Solar PV & 26/\textbf{41}/60 \\
		Wind onshore & 7/\textbf{11}/56 \\
		Wind offshore & 8/\textbf{12}/35 \\
	\end{tabular}
	\label{tab:lifecycle_emissions_per_source}
\end{table}

By calling $G$ the number of possible energy sources and $\theta_g$ the median lifecycle emission factor of source $g$, the grid's carbon intensity at time period $t$ can be calculated as:
\begin{equation}
	C_t = \sum_{g = 1}^{G} w_t^g \theta_g \quad ,
	\label{eq:grid_carbon_intensity_calculation}
\end{equation}
where $w_t^g$ is the share of energy source $g$ in the generation mix in period $t$, computed as:
\begin{equation}
	w_t^g = \frac{Q_t^g}{\sum_{h=1}^G Q_t^h} \quad ,
	\label{eq.grid_carbon_intensity_weight_factor}
\end{equation}
and $Q_t^g$ is the contribution of energy source $g$ to the grid generation mix in period $t$.
Note that this formulation expresses $C_t$ as a weighted sum of the lifecycle emission factors $\theta_g$, where the weights correspond to the shares of each energy source in the generation mix in period $t$.

\begin{figure}[H]
	\centering
	\includegraphics[width=.6\linewidth, trim={.44cm .2cm .44cm 1.2cm}, clip]{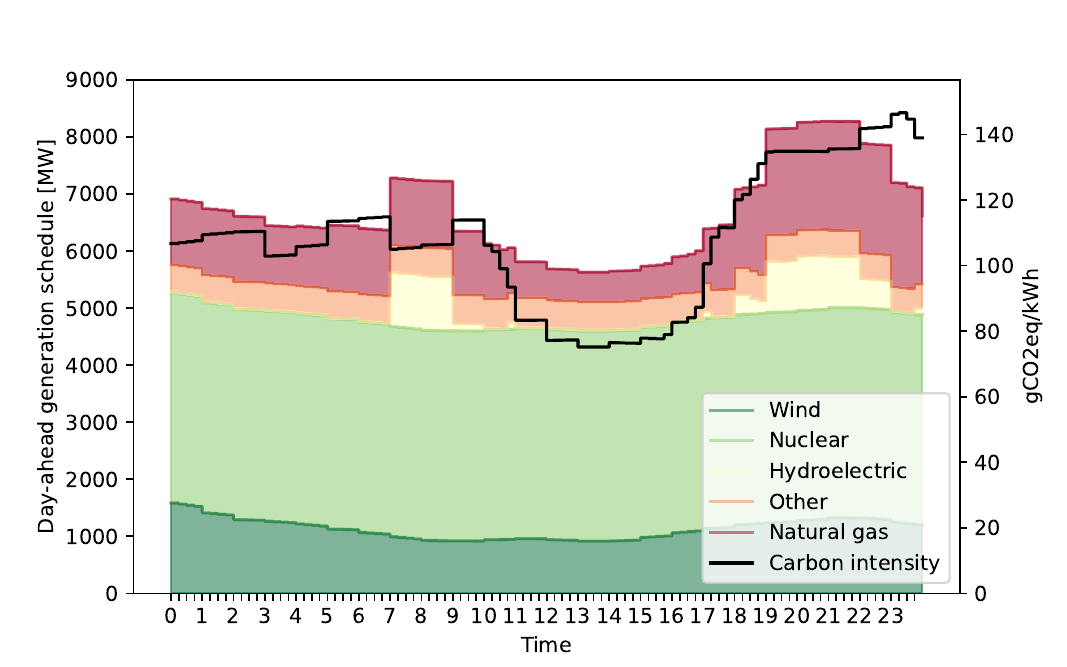}
	\caption{Grid carbon intensity and power generation mix}
	\label{fig:grid_carbon_intensity}
\end{figure}

Although the actual generation mix shifts throughout the day due to real-time adjustments for unexpected fluctuations in supply and demand, system operators establish day-ahead generation schedules that account for economic and technical constraints for plant operations.
Figure \ref{fig:grid_carbon_intensity} illustrates the day-ahead generation schedule for a representative day in Belgium, showing the planned generation mix and corresponding carbon intensity.
The carbon intensity, shown on the secondary y-axis, fluctuates significantly throughout the day due to varying generation sources.
For instance, during peak demand periods, a higher reliance on natural gas often results in increased carbon intensity.

\subsection{Carbon-aware scheduling model}
\label{subsection:CAS_model}

After discussing the time-dependent nature of job power requirements, the possible energy sources, and the methodology for computing grid carbon intensity, we now formally introduce the carbon-aware permutation flow-shop scheduling model.

Consider a set of $N$ jobs that require $M$ consecutive processing steps, each performed by a dedicated machine, over a planning horizon divided into $T$ equal time periods. The sequence of jobs is identical across all machines, following the permutation flow-shop constraint. Each operation $m \in \{1,\ldots,M\}$ of each job $i \in \{1,\ldots,N\}$ has a specific processing time of $D_{im}$ periods and requires a specific amount of electrical power $P_{im}^k$ per period, with $k \in \{1,\ldots,D_{im}\}$.

This electrical power can be sourced either from on-site renewable energy or from the public grid.
Electricity generated from on-site renewable sources is free of scope 2 GHG emissions, but its availability is both limited and time-dependent.
While on-site generation is inherently uncertain due to weather variability, planning typically relies on forecasted availability.
Based on these forecasts, the available renewable electricity in period $t$, denoted by $A_t$, is assumed to be deterministic and known for the entire planning horizon, i.e. for $t \in \{1,\ldots,T\}$.
Energy storage and feed-in of surplus electricity to the public grid are not permitted, meaning that any available renewable power is directly used to meet the schedule's electricity demand.
When on-site renewable power is insufficient, the additional electricity required is drawn from the public grid.

Unlike renewable power, grid electricity has no availability constraints.
However, its carbon intensity fluctuates over time based on the generation mix.
In this model, deterministic forecasts of the generation mix for the upcoming planning horizon are assumed to be available, allowing the carbon intensity $C_t$ for each period $t$ to be calculated using Equation (\ref{eq:grid_carbon_intensity_calculation}).

To facilitate the calculation of emission reductions due to the use of on-site renewables, we further define:
\begin{equation}
E_{im}^t = \sum_{k=1}^{D_{im}} P_{im}^k \, C_{t+k-1}
\label{eq:emissions_if_only_grid}
\end{equation}
as the scope 2 GHG emissions of scheduling operation $m$ of job $i$ to start in time period $t$ if only grid electricity is used.

Before moving on to the MILP formulation, we summarize the assumptions under which the model is formulated:
\begin{enumerate}[a)]
\item Job processing times and power requirements are deterministic and known for the upcoming planning horizon.
This is a justifiable assumption because manufacturing processes often have well-defined duration and power profiles, especially in controlled production environments.
Historical data or dry runs can provide reliable estimates.
\item On-site renewable power generation is deterministic and known for the upcoming planning horizon.
Renewable power generation can be accurately forecasted using weather data and historical generation patterns. While small deviations will exist, day-ahead predictions should be reliable for planning purposes. 
\item Grid generation mix is deterministic and known for the upcoming planning horizon.
This assumption is reasonable because generation schedules are typically planned in advance to coordinate the operation of power plants ensuring enough time to start and stop generation to meet the expected demand.
These schedules are used to provide accurate forecasts of the grid's carbon intensity.
\item On-site renewable power cannot be stored or injected into the public grid.
Therefore, whenever available, it will be used first.
This assumption aligns with scenarios where storage infrastructure is unavailable or uneconomical, and is considered for the sake of simplicity.
\end{enumerate}

To formulate the mathematical model, we define the decision variables. The main variables are the binary $x_{im}^t$, which equal $1$ if operation $m$ of job $i$ is scheduled to start at the beginning of period $t$, and $0$ otherwise. Because the different operations of a job have to be performed in the given sequence, a limited interval of periods in which operation $m$ of job $i$ can feasibly be started can be determined in preprocessing. The first feasible starting period, denoted $S_{im}$, depends on the durations of preceding operations: $S_{im} = 1 + \sum_{l=1}^{m-1} D_{il}$. Similarly, the final feasible starting period, denoted $F_{im}$, depends on the durations of remaining operations: $F_{im} = T + 1 -\sum_{l=m}^{M} D_{il}$.

To determine the job sequence and avoid overlap, the binary variables $s_{ij}$ are introduced, which equal $1$ if job $i$ comes before job $j$ in the sequence (with $i<j$). Finally, the non-negative variables $y_t$ track the on-site renewable power that is used in period $t$.

To improve readability of the model, we also introduce two sets of auxiliary variables that are derived from the $x_{im}^t$ variables. Variables $\tau_{im}$ define the period in which operation $m$ of job $i$ starts, and variables $p_{im}^t$ define the amount of electrical power consumed by operation $m$ of job $i$ in period $t$.

The MILP formulation is presented in Program (\ref{MILP}).
The objective function (\ref{MILP:obj_function}) minimizes the scope 2 GHG emissions of the schedule by first calculating the total emissions assuming exclusive use of grid electricity and then subtracting the saved emissions by using on-site generated renewable power.
This allows for pre-calculation of the emission matrix $E_{ij}^m$, which reduces computational overhead.
Constraints (\ref{MILP:schedule_each_operation}) ensure that all operations are scheduled within their feasible time horizon. Constraints (\ref{MILP:define_tau}) and (\ref{MILP:calculate_energy}) define the auxiliary variables $\tau_{im}$ and $p_{ij}^m$, respectively, based on the $x_{im}^t$ variables.

Constraints (\ref{MILP:precedence}) are the precedence constraints, preventing a job's next operation (on the next machine) from starting until the job's current operation is completed. Constraints (\ref{MILP:calculate_s_1}) and (\ref{MILP:calculate_s_2}) ensure the same sequence of jobs is maintained across all the machines using the $s_{ij}$ variables.

Finally, constraints (\ref{MILP:limit_use_of_on-site_1}) and (\ref{MILP:limit_use_of_on-site_2}) determine the total on-site renewable power used by the schedule in each period $t$, ensuring it does not exceed the available power $A_t$.

\begin{small}
\begin{mini!}|s|[2]
{\ }{\sum_{i=1}^{N} \sum_{m=1}^{M} \sum_{t=S_{im}}^{F_{im}} E^t_{im} \, x^t_{im} - \sum_{t=1}^T C_t \, y_t \label{MILP:obj_function}}{\label{MILP}}{}
\addConstraint{\sum_{t=S_{im}}^{F_{im}} x^t_{im}}{= 1 \label{MILP:schedule_each_operation}}{\enspace \begin{cases} \forall i \in \{1, \ldots, N\} \\ \forall m \in \{1, \ldots, M\} \end{cases}}
\addConstraint{\tau_{im}}{= \sum_{t=S_{im}}^{F_{im}} t \, x^t_{im} \label{MILP:define_tau}}{\enspace \begin{cases} \forall i \in \{1, \ldots, N\} \\ \forall m \in \{1, \ldots, M\} \end{cases}}
\addConstraint{p^t_{im}}{= \sum_{k=1}^{\min\left(D_{im}, t+1-S_{im}\right)} P^k_{im} \, x^{t-k+1}_{im}}{\enspace \begin{cases} \forall i \in \{1, \ldots, N\} \\ \forall m \in \{1, \ldots, M\}\\ \forall t \in \{S_{im},\ldots,F_{im}+D_{im}\} \end{cases} \label{MILP:calculate_energy}}
\addConstraint{\tau_{i(m+1)}}{\geq \tau_{im} + D_{im}}{\enspace \begin{cases} \forall i \in \{1, \ldots, N\} \\ \forall m \in \{1, \ldots, M-1\} \end{cases} \label{MILP:precedence}}
\addConstraint{\tau_{im}}{\geq \tau_{jm} + D_{jm} - T \, s_{ij}}{\enspace \begin{cases} \forall i \in \{1, \ldots, N-1\} \\ \forall j \in \{i+1, \ldots, N\} \\ \forall m \in \{1, \ldots, M\} \end{cases} \label{MILP:calculate_s_1}}
\addConstraint{\tau_{jm}}{\geq \tau_{im} + D_{im} - T \, \left( 1 - s_{ij} \right)}{\enspace \begin{cases} \forall i \in \{1, \ldots, N-1\} \\ \forall j \in \{i+1, \ldots, N\}  \\ \forall m \in \{1, \ldots, M\}\end{cases} \label{MILP:calculate_s_2}}
\addConstraint{y_t}{\leq \sum_{i=1}^{N} \sum_{m=1}^{M} p^t_{im}}{\enspace \forall t \in \{1,\ldots,T\} \label{MILP:limit_use_of_on-site_1}}
\addConstraint{y_t}{\leq A_t}{\enspace \forall t \in \{1,\ldots,T\} \label{MILP:limit_use_of_on-site_2}}
\addConstraint{x^t_{im}}{\in \{0,1\},}{\enspace \begin{cases} \forall i \in \{1, \ldots, N\} \\ \forall m \in \{1, \ldots, M\}\\ \forall t \in \{S_{im},\ldots,F_{im}\} \end{cases}}
\addConstraint{p^t_{im}}{\in \mathbb{R}^+,}{\enspace \begin{cases} \forall i \in \{1, \ldots, N\} \\ \forall m \in \{1, \ldots, M\}\\ \forall t \in \{S_{im},\ldots,F_{im}\} \end{cases}}
\addConstraint{s_{ij}}{\in \{0,1\},}{\enspace \begin{cases} \forall i \in \{1, \ldots, N-1\} \\ \forall j \in \{i+1, \ldots, N\} \end{cases}}
\addConstraint{y_t}{\in \mathbb{R}^+,}{\enspace \forall t \in \{1,\ldots,T\}} \enspace .
\end{mini!}
\end{small}

This problem can be solved with exact methods by using commercial solvers such as Gurobi or CPLEX.
However, this is only feasible for limited-sized instances.
Since an easier variant of this problem, which does not consider limited on-site electricity generation, is already NP-hard for $M \geq 3$ \cite{Garey1976}, it follows that the proposed problem is also NP-hard, at least for $M \geq 3$.
Therefore, efficient solution methods are needed to tackle real-world sized instances promptly.
In this paper, we propose a dedicated memetic algorithm, which will be discussed in the next section.

\section{Memetic algorithm framework}
\label{section:memetic_algorithm_framework}

We now introduce the dual random-key memetic algorithm designed to efficiently solve the carbon-aware permutation flow-shop scheduling problem (MA-CAS-PFSP). 
Memetic algorithms are metaheuristics that integrate local search operators within evolutionary computing, leveraging the strengths of each approach.

Evolutionary algorithms, inspired by natural evolution in biology, operate on a population of individuals, which are encoded solutions.
Each individual consists of a set of genes, where each gene corresponds to a specific property of the solution.
Through evolutionary operations such as crossover, mutation, and selection, new individuals are generated by modifying genetic material inherited from their predecessors.
Over successive generations, this process improves the overall fitness of the population.
Due to their strong global search capabilities and general applicability, evolutionary algorithms have been successfully applied to a wide range of combinatorial optimization problems, including scheduling.

However, as stated by the No Free Lunch Theorem, no single optimization algorithm is universally superior across all problems \cite{Wolpert1997}.
As a result, evolutionary algorithms can exhibit slow convergence towards optimal solutions.
In contrast, local search methods refine solutions through greedy, fine-tuning mechanisms but are prone to getting trapped in local optima \cite{Cheng1999}.
By integrating these complementary strategies, memetic algorithms have demonstrated superior performance in solving complex scheduling problems \cite{eiben2015introduction}.

\subsection{Solution representation}
The first step in evolutionary computing is defining a solution representation, which determines what information is contained in an individual and what solution properties correspond to its genes.
The solution representation serves as a link between the problem's real-world context and the algorithm's abstract problem-solving space.
The choice of data structure used for the solution representation is a critical design aspect of evolutionary algorithms, as it can significantly impact the algorithm's performance and search efficiency.
For scheduling problems, solution representations typically focus on encoding only the sequence of operations to be processed on each machine, with the assumption that operations are scheduled to start as soon as possible.
In the case of permutation flow-shop scheduling, encoding the sequence of jobs suffices, as operations are processed in the same order on all machines.

However, in carbon-aware scheduling, the solution representation must allow for planned idle times between operations, as delaying the start of an energy-intense task to a period with high on-site power generation can potentially reduce overall GHG emissions.
The proposed solution representation therefore consists of two components: 1) a job sequence, representing the order in which operations are processed on all machines, and 2) idle times, specifying planned time units of delay between the completion of one operation and the start of the next.

The importance of this dual representation lies in its ability to decouple the optimization of job sequences from the optimization of idle times.
A given job sequence may in fact result in vastly different solutions depending on how idle times are allocated.
By explicitly modeling both components, the algorithm can search for high-quality solutions more effectively by adjusting the job sequence and intermediate idle times independently.

\begin{figure}[h]
	\begin{subfigure}{.49\columnwidth}
		\centering
		\includegraphics[clip, trim={0.4cm 0cm 0.4cm 0.5cm}, width=\textwidth]{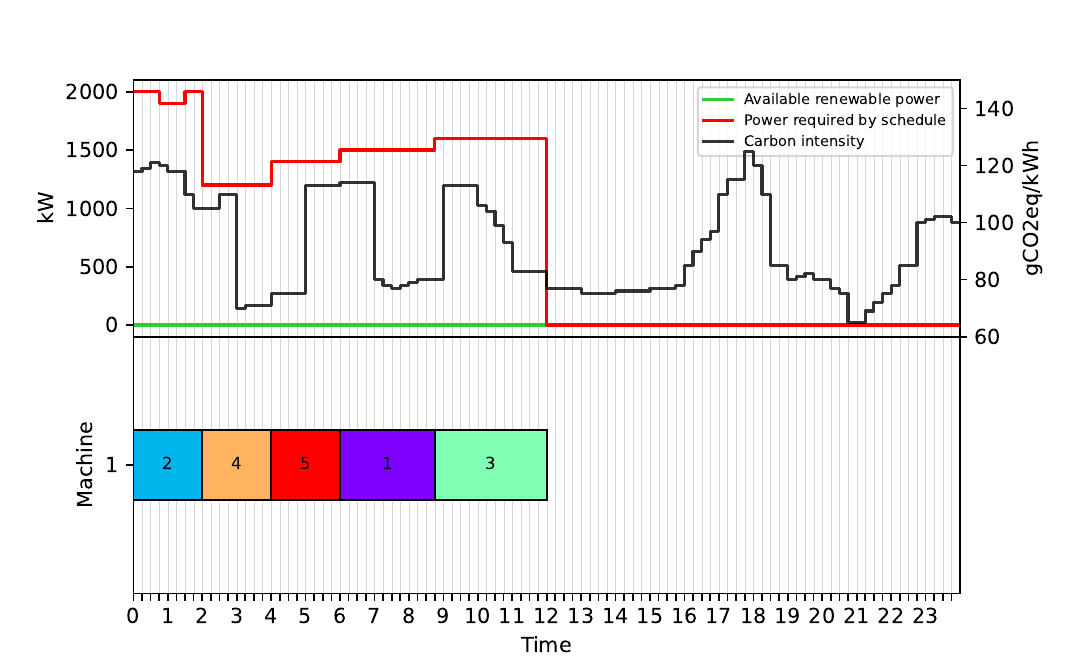}
		\caption{Jobs are scheduled as soon as possible}
		\label{fig:example_pauses_1}
	\end{subfigure}
	\begin{subfigure}{.49\columnwidth}
		\centering
		\includegraphics[clip, trim={0.4cm 0cm 0.4cm 0.5cm}, width=\textwidth]{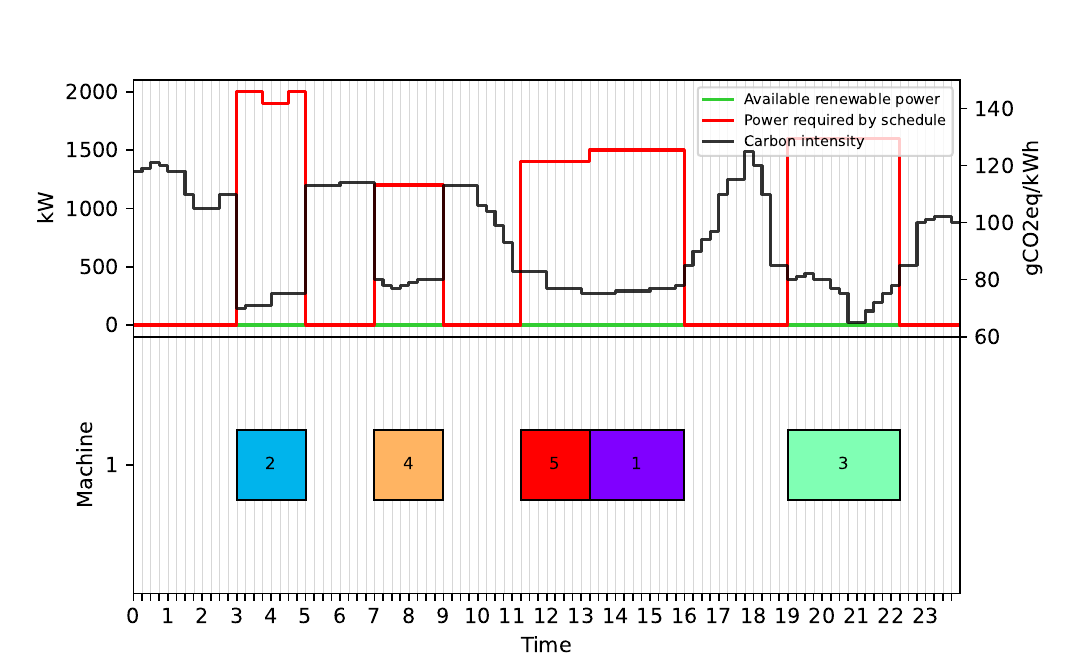}
		\caption{Pauses are incorporated between jobs}
		\label{fig:example_pauses_2}
	\end{subfigure}
	\caption{Single-machine example}
	\label{fig:example_pauses_all}
\end{figure}

To illustrate this concept, consider an example where five jobs must be scheduled on a single machine within a time horizon of $24$ hours.
Since each job consists of a single operation, the terms \textit{job} and \textit{operation} can be used interchangeably in this context.
The duration and power requirements of the jobs are known with a granularity of $15$ minutes and are shown in Table \ref{tab:example_pauses}.
\begin{table}[h]
	\fontsize{9pt}{14pt}\selectfont
	\centering
    \caption{Jobs processing times and power requirements for a single-machine example}
	\begin{tabular}{c|c|l}
		$i$ & $D_i^1$ & $\textbf{P}_i^1$ \\
		\hline
		$1$ & $11$ & $[1500,1500,1500,1500,1500,1500,1500,1500,1500,1500,1500]$  \\
		$2$ & $8$ & $[2000,2000,2000,1900,1900,1900,2000,2000]$ \\
		$3$ & $13$ & $[1600,1600,1600,1600,1600,1600,1600,1600,1600,1600,1600,1600,1600]$\\
		$4$ & $8$ & $[1200,1200,1200,1200,1200,1200,1200,1200]$  \\
		$5$ & $8$ & $[1400,1400,1400,1400,1400,1400,1400,1400]$ \\
	\end{tabular}
	\label{tab:example_pauses}
\end{table}
Figure \ref{fig:example_pauses_1} shows the Gantt chart and the corresponding power requirement profile over time for a 2-4-5-1-3 job sequence, where jobs are processed as soon as possible.
For simplicity, no on-site renewable generation is considered in this example.
This schedule yields a total GHG emission of $7.27$ tCO\textsubscript{2}.
While all the jobs are indeed completed on time, a considerable amount of slack time is present at the end of the schedule.
If this slack time is instead redistributed as planned pauses between the jobs, GHG emissions can be reduced to $5.61$ tCO\textsubscript{2} by aligning the power consumption of the schedule with periods of lower grid carbon intensity.
The total duration of the jobs is $48$ time periods of $15$ minutes each, which given the time horizon of $96$ periods results in $48$ periods of slack time. 
Figure \ref{fig:example_pauses_2} illustrates how the same job sequence can lead to lower GHG emissions by incorporating the 12-8-9-0-12-7 sequence of pauses (expressed in time periods) between jobs.

\subsubsection{Encoding and decoding strategy}
After defining the solution representation, an encoding and decoding strategy is needed to structure the information in a way that facilitates algorithmic operations.
To encode the solution representation, we employ the Random Keys (RK) method \cite{Bean1994}, where each gene corresponds to a floating-point positive number that indicates scheduling priority between jobs.
RKs have been successfully applied to evolutionary algorithms across a wide range of combinatorial optimization problems, particularly in scheduling \cite{Soares2020, Londe2025}.
A key advantage of RKs is that evolutionary operations such as crossover and mutation between individuals always produce feasible schedules for the first machine, eliminating the need for feasibility checks or repair operators.
For the remaining machines, any potential infeasibility can be easily handled in the decoding process without requiring additional repair mechanisms.
This benefit will become clearer after the decoding strategy is discussed.

As mentioned earlier, the adopted solution representation consists of two components: the job sequence and the idle times between operations on each machine.
Therefore, an encoded solution comprises two parts, one for each component.

The first part corresponds to the job sequence and consists of one array of $N$ RKs, where $N$ is the total number of jobs.
Each RK in this array represents a specific job, with the array elements ordered by the job's index $i \in \{1,\ldots,N\}$.
These RKs are normalized to sum to one.
Normalization mitigates redundancy in solution representation (i.e., multiple encodings representing the same solution), enhancing the algorithm's performance.

The second part of the solution representation consists of $M$ additional arrays, containing information about the pauses between operations on each machine.
Each of the $M$ arrays contain $(N + 1)$ RKs, which encode the lengths of pauses before, between, and after the $N$ operations, and is also normalized to one.

The decoding process transforms the RK-based representation into a feasible schedule.
This process is performed separately for each of the two parts of the solution representation.
First, the job sequence is determined by sorting the job RKs in ascending order.
The indices of the jobs, arranged by the increasing value of their corresponding RKs, define the order in which the jobs are scheduled.
Next, the idle times between operations are derived by scaling and rounding the pause RKs to match the total slack time of the schedule on machine $m$, which is precomputed as:
\begin{equation}
	S^m = T - \sum_{i=1}^{N} D_i^m \quad ,
	\label{eq:slack_calculation}
\end{equation}
where $T$ is the total available time and $D_i^m$ is the processing time of operation $O^m_i$.
The RKs in the pause arrays primarily determine the proportional distribution of slack time between consecutive operations.
Since idle times must be expressed as integers, a sum-safe rounding strategy is applied during the scaling process.
This ensures that the total assigned idle time matches the available slack time $S^m$ exactly, preventing any loss or excess due to rounding errors.
For the first machine, idle times are assigned directly.
For the remaining machines, precedence constraints must be satisfied, so the idle times are adjusted when necessary to ensure each operation starts no earlier than required by its predecessors.

\begin{figure}[H]
	\centering
	\includegraphics[width=.6\textwidth]{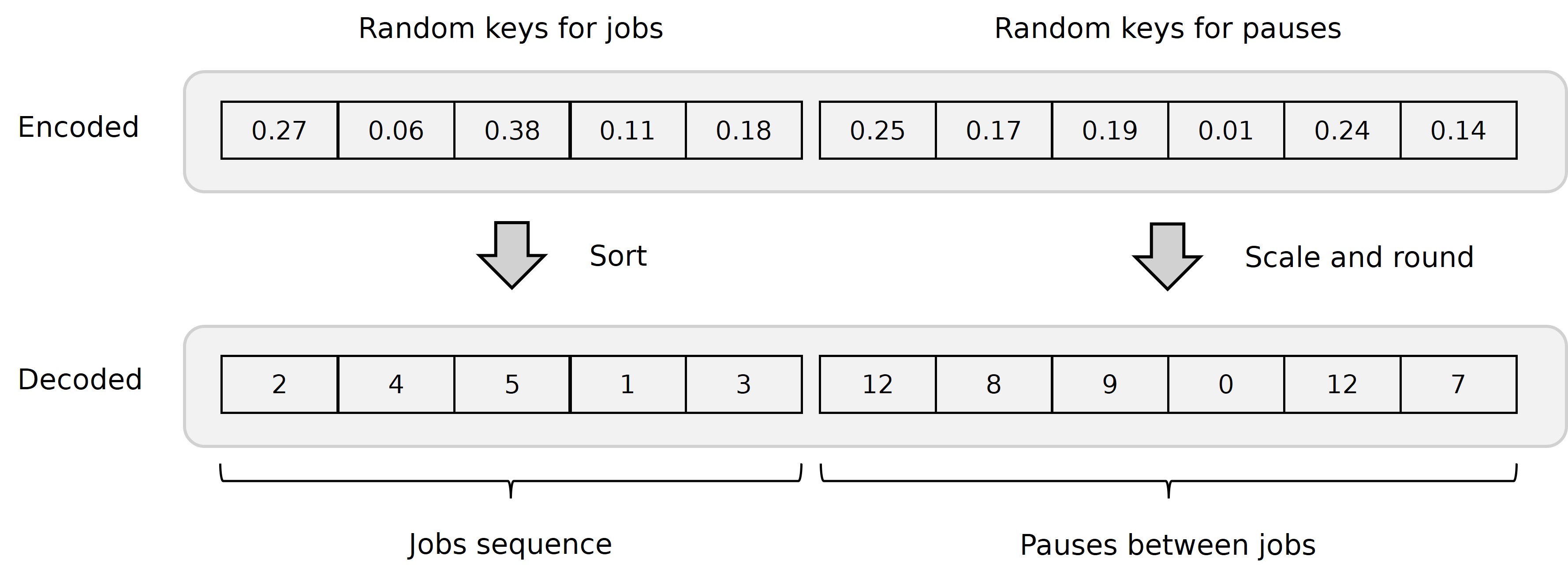}
	\caption{Encoded and decoded individual for the single-machine example of Figure \ref{fig:example_pauses_2}}
	\label{fig:example_encoding_and_decoding}
\end{figure}

To visualize the decoding process, considering again the single-machine example shown in Figure \ref{fig:example_pauses_2}, where each job consists of only one operation.
A possible encoded representation for this example and its corresponding decoded solution are depicted in Figure \ref{fig:example_encoding_and_decoding}.
As shown, the job key array encodes the priority of each job for scheduling.
Job $2$, with the smallest key value of $0.06$, is scheduled first, while job $3$, with the highest key value, is scheduled last.
Idle times are distributed based on the pause key array, with the total slack time of $48$ periods allocated proportionally.
For instance, the first idle period, corresponding to a pause key of $0.25$, receives approximately $25\%$ of the total slack time, resulting in $12$-period delay before the start of job $2$, the first job in the sequence.
Similarly, the second idle period, occurring between the end of job $2$ and the start of job $4$, receives $17\%$ of the slack time.

In this single-machine example, precedence constraints do not apply.
For multiple machines, decoding follows the same principle but includes an additional step to enforce precedence constraints.
Here, each operation's earliest possible start time is computed, and if the assigned idle time would cause an operation to begin too soon, it is delayed accordingly.

Once all operation start times are determined, the fitness of an individual is calculated by computing the net energy requirements of the schedule at each period $j \in \{1, \ldots, T\}$ and then multiplying them by the carbon intensity factor $C_j$.

As previously mentioned, adopting RKs for solution representation simplifies evolutionary operations by inherently ensuring feasibility on the first machine.
Regardless of the individual key values, they can always be sorted to determine a valid job order or scaled and rounded to compute a feasible pause sequence.
This eliminates the need for feasibility checks or repair operators, thereby improving algorithm efficiency.
In addition, the use of RKs enables effective exploration of solutions neighborhoods by providing a structured mapping between encoded and decoded solutions.
Larger changes in key values generally result in more substantial adjustments to the decoded solution, while smaller variations produce finer modifications.
This facilitates controlled adjustments during the search process, enhancing the algorithm's ability to navigate promising regions of the solution space.

\subsection{Initialization and evaluation}
The first phase of the algorithm involves creating the first population, consisting of $\rho$ individuals, and calculating their fitness level.
First generation's individuals are randomly created by sampling keys for the two types of solution arrays.
Job sequence keys are independently sampled from a uniform distribution over the interval $\left[0,1\right]$.
Once the array is complete, its keys are normalized to sum up to $1$.
Pause keys, in contrast, are sampled from an exponential distribution with expectation $1$.

The choice of an exponential distribution for pause keys is motivated by the need to increase the probability of larger values, resulting in a more asymmetric allocation of idle times.
Each machine $m \in \{1, \ldots, M\}$ has a dedicated pause key array, which, after sampling, is normalized and scaled to the available slack time on that machine.
This approach results in a more varied spread of idle times between operations compared to a uniform distribution.
In particular, it increases the likelihood of generating solutions where idle time intervals occasionally span multiple consecutive periods, which better aligns with the characteristics of energy-related data, as periods of low on-site renewable generation or high grid carbon intensity often extend over multiple time units.

While the use of RKs guarantees feasibility of the schedule on the first machine, this does not necessarily hold for subsequent machines.
Due to precedence constraints between consecutive operations of the same job, a schedule that is feasible on the first machine may cause the time limit $T$ to be exceeded on subsequent machines, making the schedule unfeasible.
This is particularly likely when a considerable idle time is planned at the beginning of the schedule.
Instead of repairing such infeasible schedules, which would require additional computations, we assign them a penalty proportional to their lateness.
Specifically, for a given schedule, we compute the completion time $CT$ of the last operation on the last machine and apply the following penalty:
\begin{equation}
	\text{Penalty} = \max\left(0, \; CT - T\right) \times 10^{10} \enspace .
\end{equation}
This penalty discourages infeasible schedules while allowing the evolutionary process to naturally eliminate such low-fitness individuals over successive generations.

To ensure that at least one feasible schedule is present in the initial population, we first generate $\rho - 1$ random individuals and then manually include a single trivial feasible schedule where jobs are processed in a first-come-first-serve (FCFS) sequence with no planned idle time between operations.
With the adopted solution representation, an FCFS schedule is naturally obtained by using a job sequence array of $N$ ascending RKs, and $M$ pause arrays where only the $(N + 1)$-th element is nonzero.
With this structure, our decoding strategy effectively assigns all slack time after the completion of the last operation on each machine.
Combined with the precedence constraints, this ensures that the resulting schedule follows the FCFS principle.
While not strictly necessary, this heuristic initialization accelerates convergence with only minimal additional computational overhead.
The initialization procedure is then finalized by decoding the individuals and calculating their fitness, resulting in the first current generation of the algorithm.

\subsection{Controlled swap crossover}
Once the current generation is available -- whether from the algorithm's initialization or the start of a new iteration -- the algorithm begins generating offspring.
The first operation in this process is crossover, where new individuals are generated by recombining the genetic information of parent individuals.

We propose a controlled swap crossover operator with uniform parent selection.
Each individual in the current generation has an equal probability of being selected as a parent, with the condition that the two chosen parents must be distinct.
Once a pair of parents is chosen, a subset of their corresponding keys are swapped to produce two new individuals.
This swapping process is performed independently for the job sequence key array and the pause key arrays, and is guided by two types of parameters that control both the frequency and the intensity of the crossover.

The first type is the crossover rate $\xi$, which determines the fraction of offspring generated through crossover.
The remaining $(1 - \xi)\%$ of individuals are selected directly from the current generation using an elitist approach.
For example, with a population size of $\rho=100$ individuals and $\xi=0.7$, $70$ individuals in the offspring will be generated through crossover.

The second type is the crossover probability $\chi$, which defines the probability of swapping each gene between parents within a given array.
This controls the severity of the crossover.
For instance, with $N=5$ jobs and $\chi=0.4$, each job key has a $40\%$ chance of being swapped, meaning that, on average, $2$ out of $5$ job keys will be exchanged per crossover operation.

To account for the differing roles of the job sequence array and pause arrays, separate crossover probabilities, $\chi_j$ and $\chi_p$, are defined.
These parameters independently control the swapping likelihood for the job sequence and the pause keys, respectively.

\begin{figure}[H]
	\centering
	\includegraphics[width=.6\textwidth]{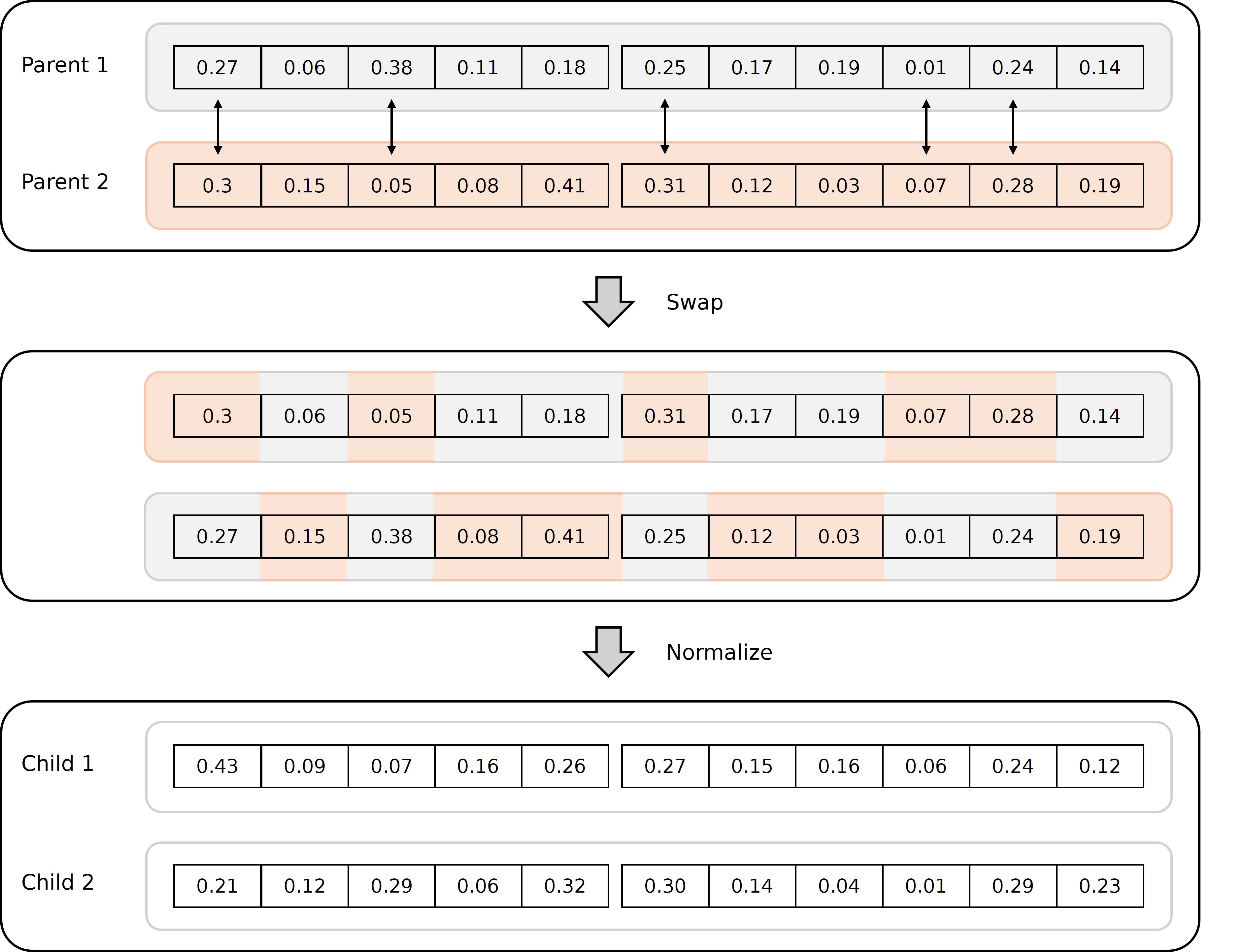}
	\caption{Crossover on a single-machine example with $N=5$, $\chi_j=0.4$, and $\chi_p=0.5$}
	\label{fig:crossover}
\end{figure}

After swapping, the key arrays are normalized to sum up to one.
The generation of two new individuals with the proposed crossover operator is illustrated in Figure \ref{fig:crossover}.

\subsection{Nonuniform mutation}
The second operation applied to the newly created offspring is mutation.
We employ a nonuniform mutation strategy, where the genes of an individual are modified by adding a random value to its RKs.
These random values are independently sampled from a normal distribution with mean of zero and a variable standard deviation: $\sigma_j$ for the job key array and $\sigma_p$ for the pause key arrays.
These standard deviations, referred to as mutation step sizes, control the severity of the mutation, determining how far the new value deviates from the original one.
Smaller step sizes produce mutations closer to the original value, while larger step sizes introduce more significant changes.
To ensure all keys remain non-negative after mutation, any resulting negative values are clipped to zero.
The key arrays are then again normalized, preserving the consistency of the representation.

Although every individual in the offspring undergoes the mutation process, not all the genes within an individual are necessarily altered.
For each individual, the mutation is applied probabilistically to each gene, with the likelihood of mutation determined by the respective mutation probabilities: $\pi_j$ for the job keys and $\pi_p$, for the pause keys.
This means that each gene in an individual is independently evaluated against its mutation probability to decide whether it will be mutated.
With these four parameters, the frequency ($\pi_j$ and $\pi_p$) and severity ($\sigma_j$ and $\sigma_p$) of mutation for job sequence keys and pause keys can be independently controlled.

\begin{figure}[H]
	\centering
	\includegraphics[width=.6\textwidth]{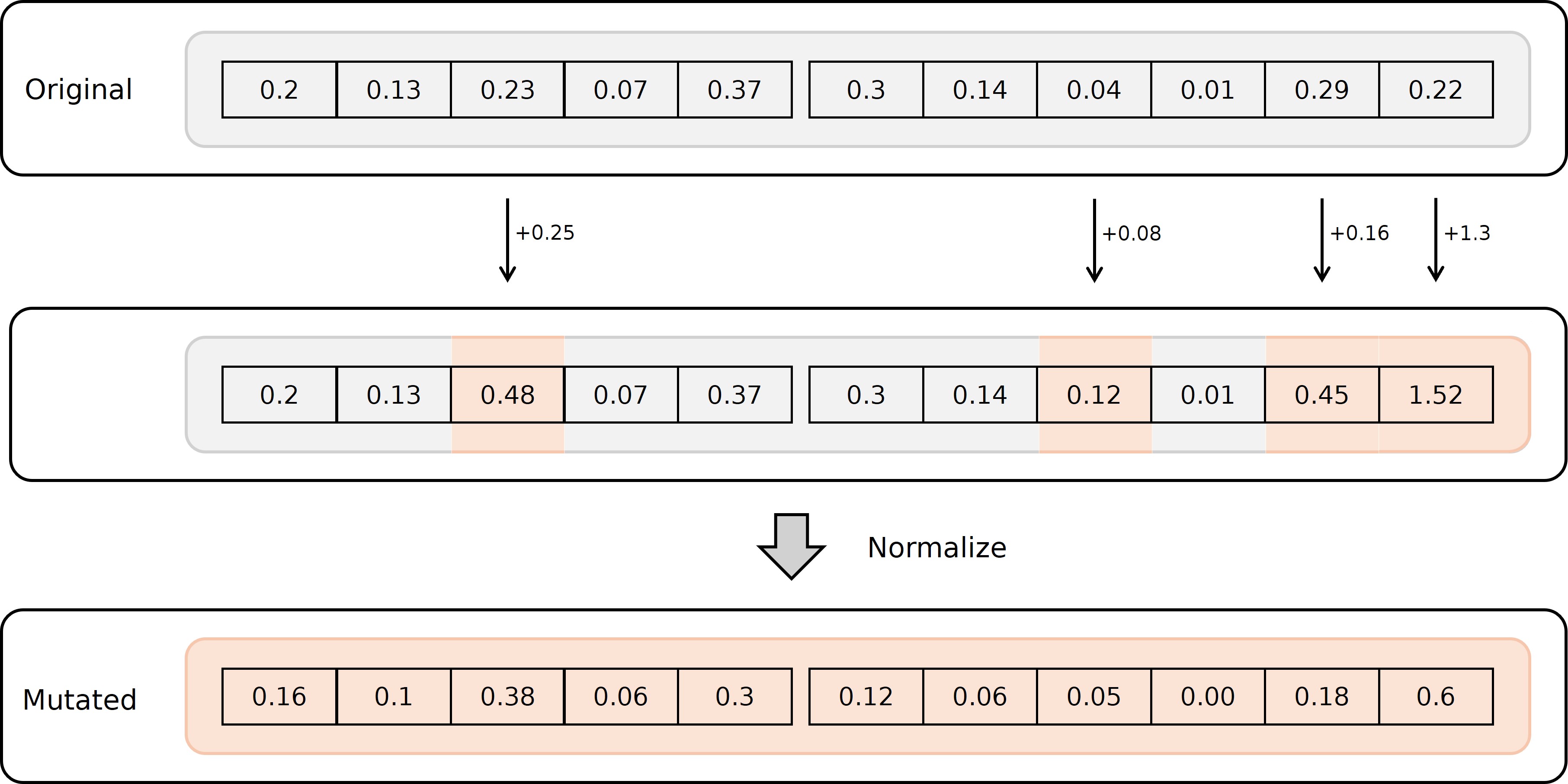}
	\caption{Mutation on a single-machine example with $N=5$, $\pi_j=0.2$, $\sigma_j=0.06$, $\sigma_p=0.2$, and $\pi_p=0.5$}
	\label{fig:mutation}
\end{figure}

The effect of mutation on a single individual is illustrated in Figure \ref{fig:mutation}, where a single-machine example with $N=5$ jobs is considered.
For the job key array, the mutation probability is $\pi_j=0.2$, meaning that each job key has a $20\%$ chance of being mutated.
As a result, on average, one out of the five job keys will undergo mutation.
The severity of the mutation is controlled by the step size $\sigma_j=0.06$, indicating that a random value sampled from $\mathcal{N}\left(0,\sigma_j\right)$ is added to the selected gene.
For the pause key array, the mutation probability is set at $\pi_p=0.5$, which means that, on average, half of the pause genes -- three out of the six in this example -- will be mutated.
The larger step size of $\sigma_p=0.2$ results in more significant changes for the mutated pause keys compared to the job keys.

\subsection{Local search}
\begin{figure}[h]
	\centering
	\includegraphics[width=.6\linewidth]{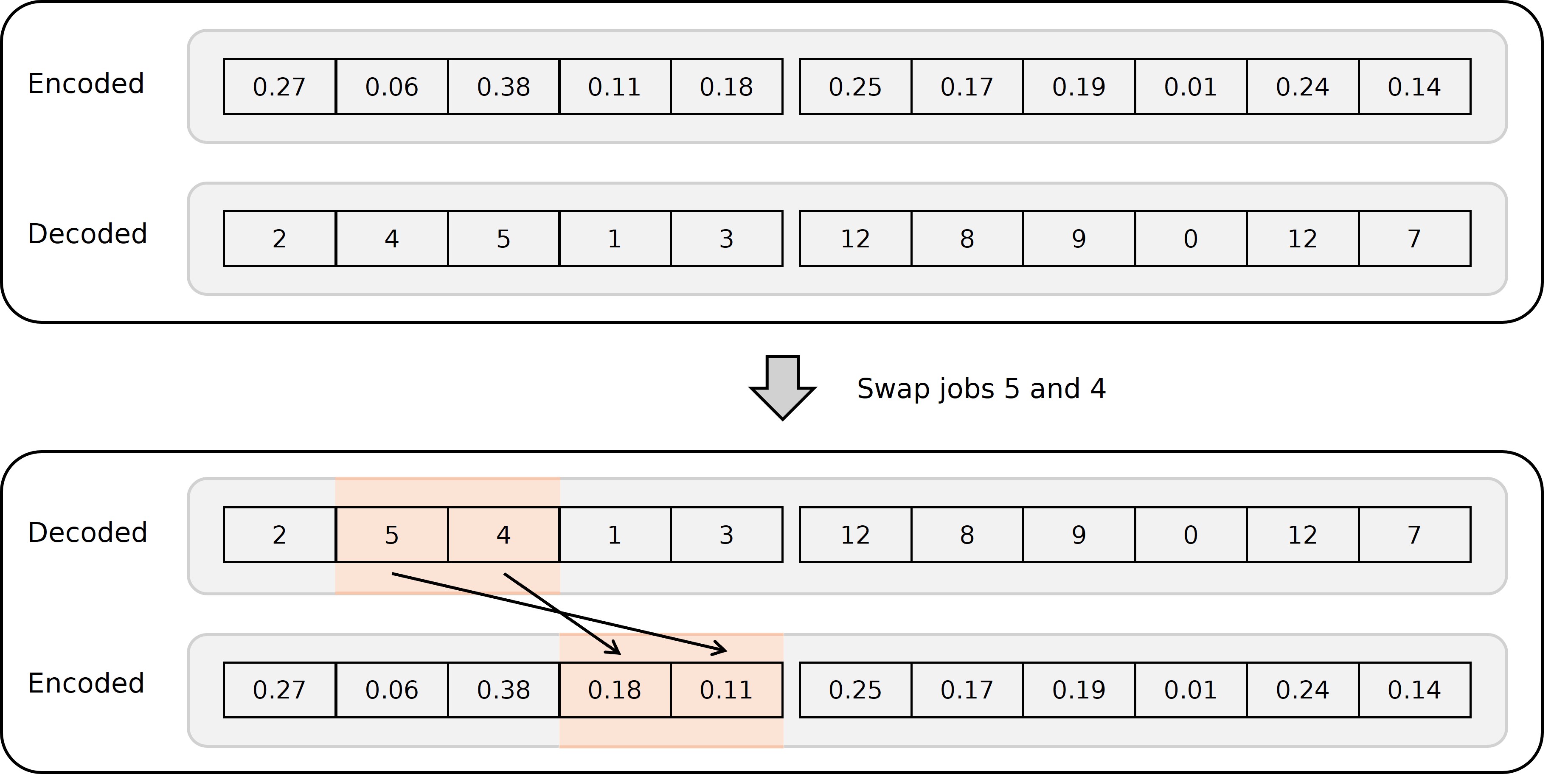}
	\caption{Local search on a single-machine example individual}
	\label{fig:overview_localsearch}
\end{figure}

Next, local search is performed to improve the solution by adapting the job sequence.
We adopt a pairwise adjacent swap operator, which iteratively examines adjacent job swaps in the sequence and recalculates the fitness after each swap.
To balance solution quality with computational efficiency, the search stops as soon as an improvement is found.
At this point, the RKs of the corresponding swap are also exchanged, resulting in an encoded improved solution.
Figure \ref{fig:overview_localsearch} shows an example of one possible swap with the proposed local search operator.

\subsection{Selection and termination criterion}
After crossover, mutation, and local search, the creation of offspring is completed.
Next, the current population and offspring are merged into a candidate pool, where individuals are ranked in descending order of fitness.
The top $\rho$ individuals are then selected to form the new generation.
This elitist selection strategy ensures that only the most fit individuals, including any improvements over the current population, are retained.
The process is repeated until the predefined maximum number of generations $\gamma$ is reached.

\begin{figure}[h]
	\centering
	\includegraphics[width=.8\linewidth]{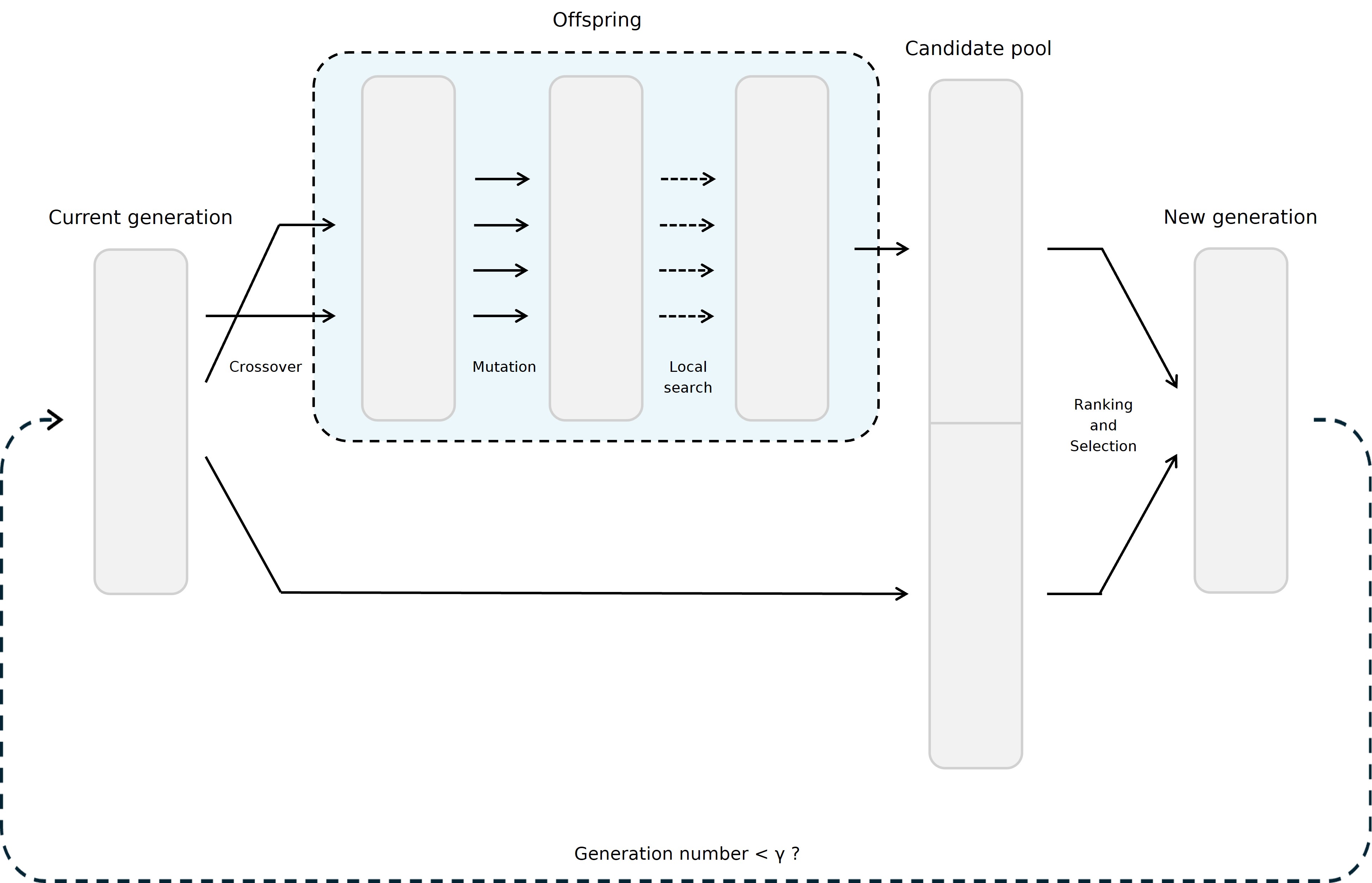}
	\caption{Overview of the proposed MA-CAS-PFSP algorithm}
	\label{fig:overview_algorithm}
\end{figure}

An overview of the proposed MA-CAS-PFSP algorithm framework is presented in Figure \ref{fig:overview_algorithm}, while Table \ref{tab:algorithm_parameters} summarizes its relevant parameters and their definition.

\begin{table}[h]
	\fontsize{9pt}{14pt}\selectfont
	\centering
    \caption{Summary algorithm parameters and their definitions}
	\begin{tabular}{c|l}
		Parameter & Definition\\
		\hline
		$\rho$ & Population size  \\
		$\gamma$ & Maximum number of generations \\
		$\xi$ & Crossover rate  \\
		$\chi_j$ & Crossover probability for jobs  \\
		$\chi_p$ & Crossover probability for pauses \\
		$\pi_j$ & Mutation probability for jobs \\
		$\pi_p$ & Mutation probability for pauses \\
		$\sigma_j$ & Mutation step size for jobs \\
		$\sigma_p$ & Mutation step size for pauses \\
	\end{tabular}
	\label{tab:algorithm_parameters}
\end{table}

\section{Computational experiments}
\label{section:computational_experiments}
We now present and discuss the computational experiments conducted in this study by first detailing the instance dataset generation procedure.
Next, we outline the algorithm parameter tuning process.
Finally, we discuss the results of computational experiments for both single-machine and multi-machine models.

\subsection{Instance dataset}
\begin{figure}[H]
	\centering
	\includegraphics[width=.5\linewidth]{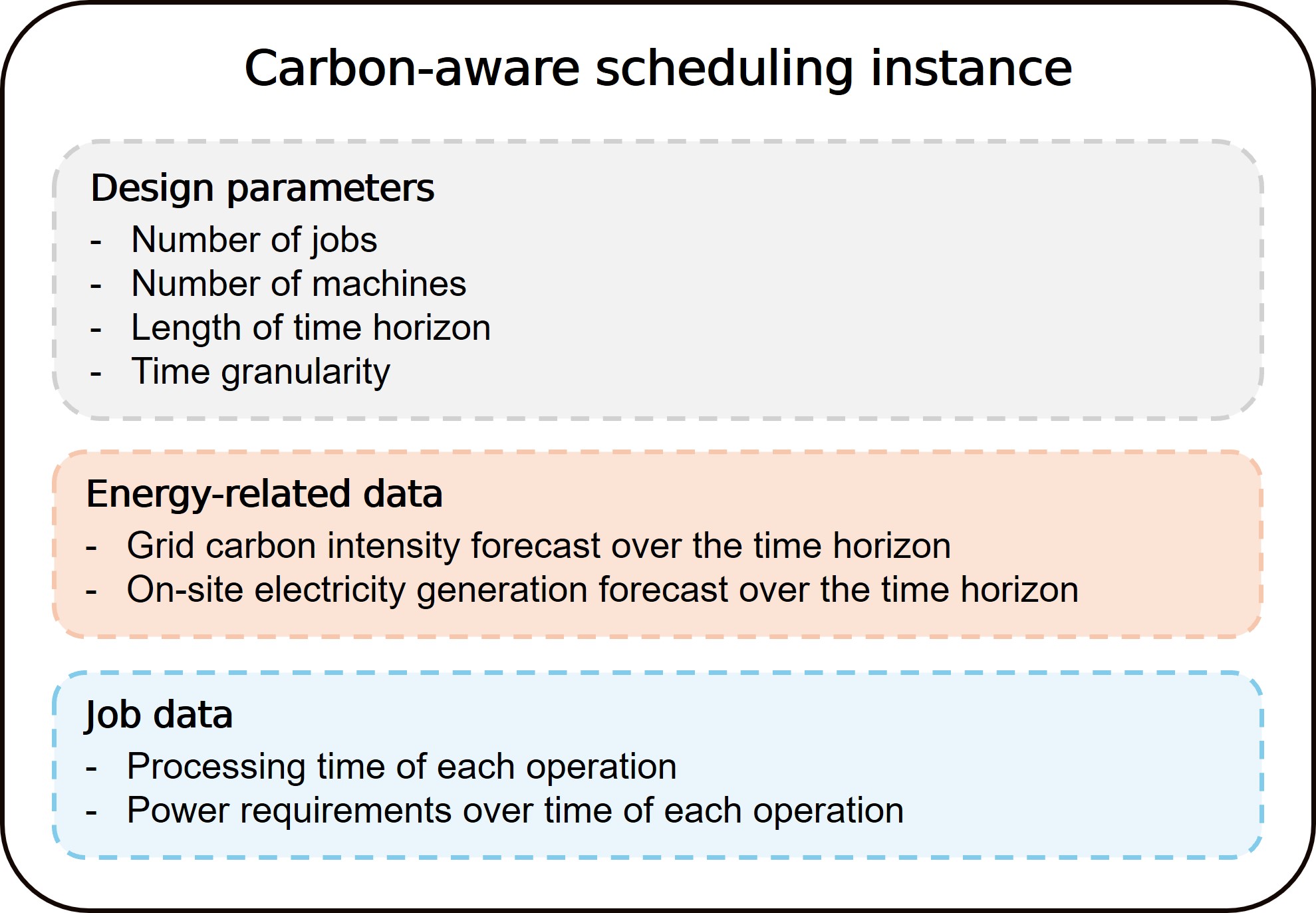}
	\caption{Overview of the information contained in each instance}
	\label{fig:instance_data}
\end{figure}

To evaluate the model and compare the algorithm's performance under varying conditions, a dataset of random problem instances is required.
As benchmarks for carbon-aware scheduling are, to the best of our knowledge, not yet available in the literature, we created a custom dataset of random instances based on historical and self-generated data.
An overview of the information stored in each instance is given in Figure \ref{fig:instance_data}.

The instance generation procedure started by defining the experiment design parameters, including the number of machines, the length of the time horizon expressed in planning days, and the time granularity.
The creation of the remaining data required for each instance is discussed in the following sections.

\subsubsection{Grid carbon intensity forecasts}
\begin{figure}[H]
	\centering
	\includegraphics[width=.8\linewidth, trim={6.5cm 5.5cm 4.5cm 5.5cm}, clip]{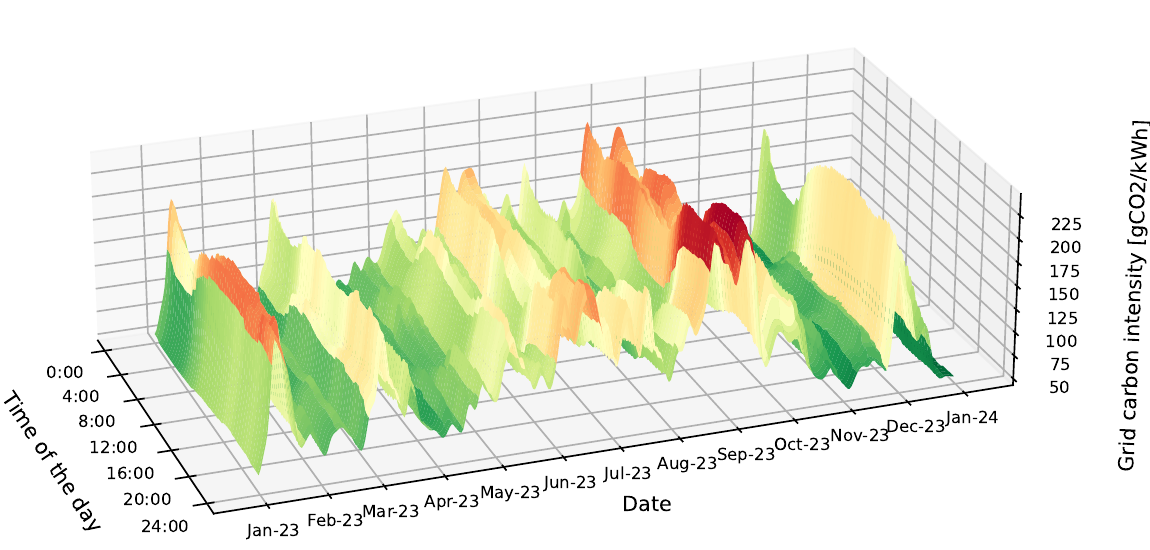} 
	\caption{Hourly carbon intensity of the Belgian grid in 2023}
	\label{fig:carbon_intensity_3dplot}
\end{figure}

Grid carbon intensity forecasts were generated using one year of historical generation mix data for the Belgian electricity grid.
This dataset, obtained from the Open Data Platform of Elia \cite{Elia_DataSet_GridMix}, the Belgian transmission system operator, provides data at a 15-minute resolution, which we used as the time granularity for all instances.
The grid's carbon intensity over time was estimated using Equation (\ref{eq:grid_carbon_intensity_calculation}), applying median lifecycle emission factors for each power source as detailed in Section \ref{subsection:electricity_demand_and_grid_carbon_intensity}.
The one-year dataset was divided into individual days, allowing random selection of a planning horizon for each instance.
The selected data served as the carbon intensity forecast for the instance's planning horizon.
An example of grid carbon intensity profile for one specific day was already shown in Figure \ref{fig:grid_carbon_intensity}, while Figure \ref{fig:carbon_intensity_3dplot} illustrates an overview of the hourly carbon intensity of the Belgian grid in $2023$, which served as input data for this study.

\subsubsection{On-site electricity generation forecasts}
\begin{figure}[H]
	\centering
	\includegraphics[width=.8\linewidth, trim={6.5cm 5.5cm 4.5cm 5.5cm}, clip]{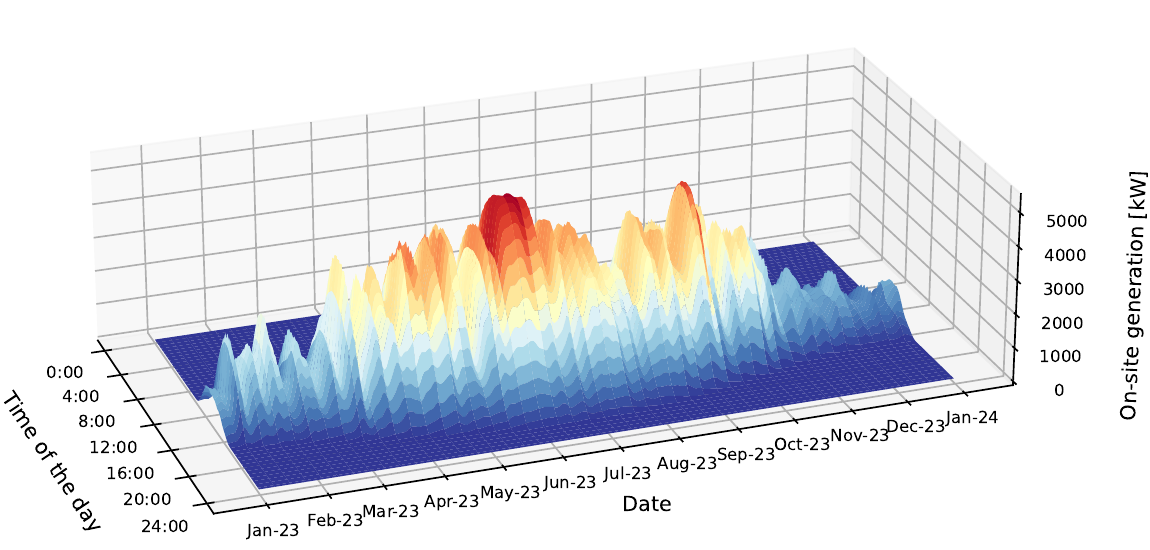} 
	\caption{Scaled hourly on-site solar generation in East Flanders, Belgium, in 2023}
	\label{fig:on-site_generation_3dplot}
\end{figure}

On-site electricity generation forecasts were also sourced from the Open Data Platform of Elia \cite{Elia_DataSet_Onsite}, focusing specifically on solar power generation.
Assuming that solar production patterns remain consistent over the years, we used one year of historical solar power generation data from East Flanders, Belgium, as the basis for the forecast.
To adapt this data for a single production site, we scaled it by $0.005$, ensuring that the magnitude of on-site generation remained comparable to job power requirements.
Although using solar generation from a broader region as a proxy for a specific site may reduce the impact of local variability, it still captures overall fluctuations in renewable availability.
Similar to the grid carbon content, this historical data was treated as forecasts for the planning horizon in each instance.
Figure \ref{fig:on-site_generation_3dplot} provides an overview of the scaled hourly on-site generation in East Flanders in $2023$, used as input data for this study.

\subsubsection{Jobs processing times and power requirements}
The generation of job data began with the creation of a pool of available operations.
As outlined in Section \ref{section:model_formulation}, each job comprises an equal number of operations $M$, with each operation having its own processing time and power requirements over time.

The processing times of operations, expressed in equal time periods, were sampled from a discrete uniform distribution.
Considering the quarter-hour granularity of the available historical data, the sampling intervals were defined as $\{2, \, 16\}$ for single-machine instances, and $\{0, \, 8\}$ for multiple-machine instances.
This approach allows for the inclusion of dummy operations with zero processing time on certain machines for multiple-machine instances, effectively representing jobs that do not require processing on those machines.

While determining the processing times for each operation, power requirements were also assigned.
First, a base power requirement was sampled from a discrete uniform distribution within the interval $\{100, \, 3000\}$.
Subsequently, for each time period of the operation, an individual power requirement was calculated by adding a random value sampled from $\mathcal{U}\{-250, \, 250\}$ to the operation's base power requirement, ensuring that the resulting requirements remained non-negative.
This approach introduces variability both between operations and across time periods within a single operation, thereby modeling realistic fluctuations in power requirements of jobs, as discussed in Section \ref{subsection:electricity_demand_and_grid_carbon_intensity}.

\subsubsection{Instance generation}
\label{subsubsection:instance_generation}
Once all required data were prepared, individual problem instances were created.
The process began by selecting a random time window from the historical energy-related data, corresponding to the planning horizon $T$.
The selected data provided the grid's carbon intensity and on-site generation forecasts for the chosen horizon.
Then, job data was added to the instance by randomly combining operations from the operations pool into jobs, making sure that the time needed to process them did not exceed the available time.
The process was repeated until four datasets of fifty random instances each were created.
An overview of the instance generation procedure is shown in Algorithm \ref{alg:instance_generation}.
Table \ref{tab:overview_instance_datasets} summarizes the instance datasets created for the experiments conducted in this study along with their key features.
In the table, $|O|$ denotes the total number of operations in an instance, while $S^{\bar{m}}$ represents the median slack time per machine.
The full instance files are available at \url{https://github.com/ugent-isye/CAS-PFSP}. 

\begin{table}[h]
	\fontsize{9pt}{14pt}\selectfont
	\centering
    \caption{Overview of the instance datasets and their key features}
	\begin{tabular}{c|c|c|c|c|c}
		Name & \# instances & $M$ & $T$ & $|O|$ (Min/Median/Max) & $S^{\bar{m}}$ (Min/Median/Max) \\
		\hline
		CAS-PFSP-M1T1 & $50$ & $1$ & 96 & 6/\textbf{10}/15 & 1/\textbf{6}/15 \\
		CAS-PFSP-M1T3 & $50$ & $1$ & 288 & 25/\textbf{32}/40 & 1/\textbf{6}/16 \\
		CAS-PFSP-M3T1 & $50$ & $3$ & 96 & 24/\textbf{36}/54 & 24/\textbf{47}/67 \\
		CAS-PFSP-M3T3 & $50$ & $3$ & 288 & 102/\textbf{135}/183 & 66/\textbf{106}/154 \\
	\end{tabular}
	\label{tab:overview_instance_datasets}
\end{table}

\begin{algorithm}
	\scriptsize
	\caption{Instance generation algorithm}\label{alg:instance_generation}
	\KwData{\(M\), \(T\), \(\mathbf{C}_{hist}\), \(\mathbf{A}_{hist}\)}
	\KwResult{\(\mathbf{I}\): instance dataset}
	\Comment{Create pool of operations}
	$o \gets 1$\; 
	$\mathbf{O} \gets$ [ ]\;
	\While{$o < 2000$}{
			\Comment{Sample processing time}
			$D_o \gets$ random draw from $\mathcal{U}\{0, \, 8\}$\; 
			\Comment{Sample power requirements}
			$P \gets$ random draw from $\mathcal{U}\{100, \, 3000\}$\;
			$k \gets 0$\;
			\While{$k < D_o$}{
				$\Delta P \gets$ random draw from $\mathcal{U}\{-250, \, 250\}$\;
				$P_{ok} \gets \max \left(0, \, P + \Delta P \right)$\;
				$k \gets k + 1$\;
				}
			Append $\left[\mathbf{P}_o\right]$ to $\mathbf{O}$\; 
			$o \gets o + 1$
		}
	\Comment{Create instances}
	$n \gets 1$\; 
	\While{$n \leq 50$}{
		$t \gets 0$\; 
		$i \gets 1$\; 
		\While{$t < T$}{
			\Comment{Add jobs to instance}
			CandidateJob $\gets$ [ ]\;
			$m \gets 1$\; 
			\While{$m \leq M$}{
				\Comment{Pick an operation from the pool}
				$o \gets$ random draw from $\mathcal{U}\{1, \, 2000\}$\;
				$O^m_i \gets \mathbf{O}_o$\;
				Append $O^m_i$ to CandidateJob\;
				$m \gets m + 1$\;
			}
			\Comment{Check feasibility}
			$t_{new} \gets$ FCFS makespan with CandidateJob\;
			\If{$t_{new} < T$}{
				$\mathbf{J}_i \gets$ CandidateJob\;
				$t \gets t_{new}$\;
				}
			$i \gets i + 1$\;
			}
		\Comment{Add energy-related data}
		$d \gets$ random draw from $\mathcal{U}\{1, \, 365\}$\; 
		$\mathbf{C} \gets \mathbf{C}_{hist}[d:d+T]$\;
		$\mathbf{A} \gets \mathbf{A}_{hist}[d:d+T]$\;
		$\mathbf{I}_n \gets \left(\mathbf{J}, \mathbf{C}, \mathbf{A}\right)$\;
		$n \gets n + 1$\;
		}
\end{algorithm}

\subsection{Parameter tuning}
Parameter tuning was performed in Python with Optuna \cite{Optuna2019}, a state-of-the-art open-source hyperparameter optimization framework.
Optuna employs Tree-structured Parzen Estimator (TPE) \cite{Bergstra2011}, a Bayesian optimization method, to model the performance of different parameter configurations and guide the search towards promising regions of the search space.
For each dataset, a separate tuning set of $10$ additional instances was generated using the same procedure, but kept apart from the test instances.
The tuning process consisted of $1000$ trials per dataset, with Optuna automatically proposing a new parameter configuration in each trial based on the outcomes of previous evaluations.
For each configuration, the algorithm was run on all instances of the tuning set, and the average objective value across these runs was returned as feedback to guide the tuning process.
Throughout the process, the population size was fixed at $\rho = 250$, and the maximum number of generations was set to $\gamma = 100$.
These values were chosen based on preliminary testing and were found to provide a good balance between solution quality and computational time.
The parameter configuration that yielded the best average performance on the tuning set was then adopted as the optimal setting for the corresponding dataset.
An overview of the results is provided in Table \ref{tab:overview_parameter_tuning}.
\begin{table}[h]
	\fontsize{9pt}{14pt}\selectfont
	\centering
    \caption{Overview of the best found parameter values for each dataset}
	\begin{tabular}{c|c|c|c|c}
		 & CAS-PFSP-M1T1 & CAS-PFSP-M1T3 & CAS-PFSP-M3T1 & CAS-PFSP-M3T3 \\
		\hline
		$\xi$ & $0.5851$ & $0.5565$ & $0.8273$ & $0.8203$  \\
		$\chi_j$ & $0.3779$ & $0.1168$ & $0.3596$ & $0.4297$  \\
		$\chi_p$ & $0.1041$ & $0.4627$ & $0.2963$ & $0.0681$  \\
		$\pi_j$ & $0.1662$ & $0.0589$ & $0.0679$ & $0.0113$  \\
		$\pi_p$ & $0.1985$ & $0.0227$ & $0.0330$ & $0.0084$  \\
		$\sigma_j$ & $0.0564$ & $0.0168$ & $0.1039$ & $0.0050$  \\
		$\sigma_p$ & $0.1873$ & $0.1832$ & $0.1959$ & $0.1901$  \\
	\end{tabular}
	\label{tab:overview_parameter_tuning}
\end{table}

\subsection{Experiments design}
\label{subsection:experiments_design}
The experiments aimed to evaluate the performance and scalability of the proposed MA-CAS-PFSP algorithm across various instance sizes.
To achieve this, we conducted four sets of experiments, each corresponding to one of the datasets discussed in Section \ref{subsubsection:instance_generation}.

In each experiment, the MA-CAS-PFSP algorithm was executed $10$ times on each of the $50$ problem instances within the dataset.
Since the search process is inherently random, multiple runs are necessary to ensure variability in the results is effectively captured.

For each problem instance, we report the average objective values achieved across the $10$ runs of the memetic algorithm, along with the average computation time.
To provide a normalized measure of variability in results, the coefficient of variation (CV) is also included, defined as the ratio of the standard deviation to the mean objective value.

The performance of the proposed algorithm is assessed by comparing the results to those generated with the CPLEX implementation of the MILP model presented in Section \ref{section:model_formulation}.
For all experiments, the population size was set to $\rho = 250$ individuals, and the algorithm was terminated after $\gamma = 100$ generations.
To ensure a fair comparison, a time limit of $1$ minute was applied to each run of both MA-CAS-PFSP and CPLEX.
This time limit was chosen to provide both methods with a similar computational budget, as MA-CAS-PFSP required less than one minute per run for all problem instances.
Additionally, to investigate the convergence speed of the exact method, we also report CPLEX results obtained with the longer time limit of $30$ minutes.

To assess the algorithm's performance, we define the percentage gap as:
\begin{equation}
	\% \Delta = \frac{\left( \text{Obj\textsubscript{CPLEX} - Obj\textsubscript{MA}} \right)}{\text{Obj\textsubscript{CPLEX}}} \cdot 100 \enspace ,
\end{equation}
where Obj\textsubscript{CPLEX} and Obj\textsubscript{MA} refer to the average objective values found within the given time limit with CPLEX and MA, respectively.

For the CPLEX results, the best feasible solution found within the specified time limit was recorded for comparison.

All the experiments were conducted on a machine equipped with an 11th Gen Intel(R) Core(TM) i5-1145G7 processor (2.60 GHz), using Python version $3.10.0.$rc2 and CPLEX version $22.1.1.0$.

\subsection{Results}
\label{subsection:results}
We now report the results of the experiments described in Section \ref{subsection:experiments_design}, highlighting general trends across the four datasets.

\subsubsection{Solution characteristics}
\label{subsubsection:solution_characteristics}
\begin{figure}[h]
	\begin{subfigure}{.49\columnwidth}
		\centering
		\includegraphics[clip, trim={0.4cm 0cm 0.4cm 0.5cm}, width=\textwidth]{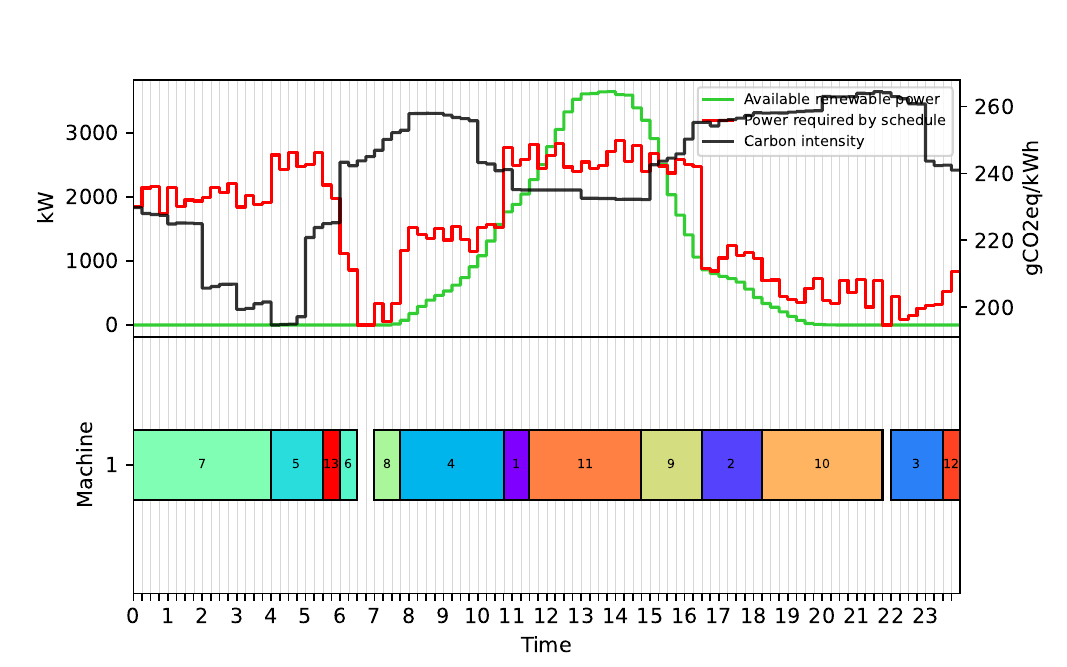}
		\caption{CAS-PFSP-M1T1}
		\label{fig:results_gantt_M1T1}
	\end{subfigure}
	\begin{subfigure}{.49\columnwidth}
		\centering
		\includegraphics[clip, trim={0.4cm 0cm 0.4cm 0.5cm}, width=\textwidth]{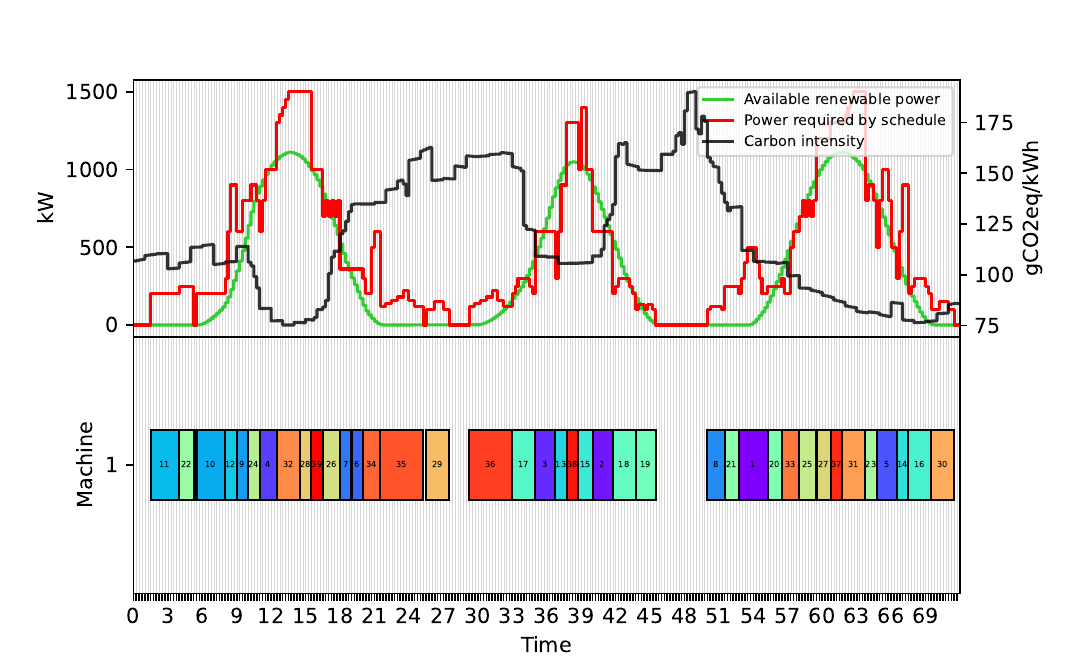}
		\caption{CAS-PFSP-M1T3}
		\label{fig:results_gantt_M1T3}
	\end{subfigure}
	\begin{subfigure}{.49\columnwidth}
		\centering
		\includegraphics[clip, trim={0.4cm 0cm 0.4cm 0.5cm}, width=\textwidth]{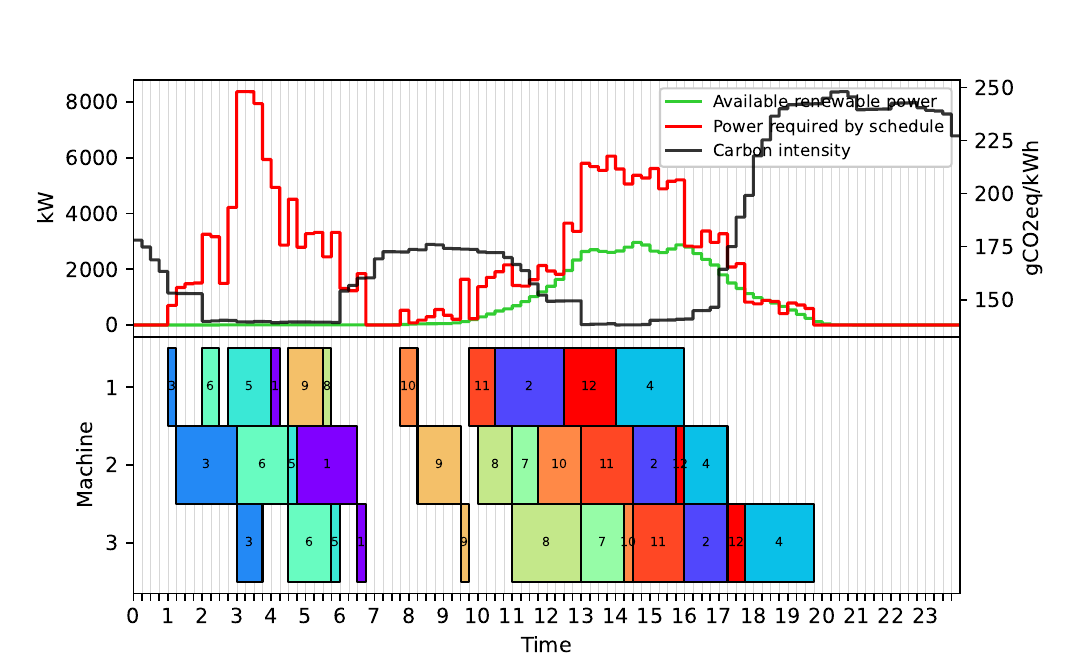}
		\caption{CAS-PFSP-M3T1}
		\label{fig:results_gantt_M3T1}
	\end{subfigure}
	\begin{subfigure}{.49\columnwidth}
		\centering
		\includegraphics[clip, trim={0.4cm 0cm 0.4cm 0.5cm}, width=\textwidth]{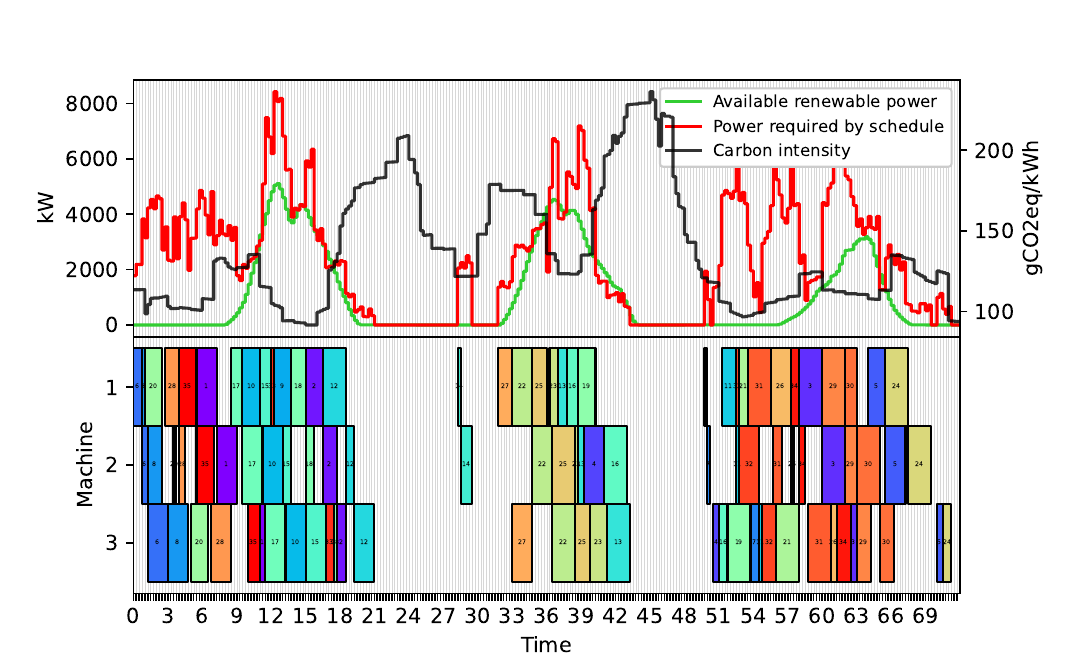}
		\caption{CAS-PFSP-M3T3}
		\label{fig:results_gantt_M3T3}
	\end{subfigure}
	\caption{Gantt chart and power requirement profiles for a representative instance of each dataset}
	\label{fig:results_gantt_all}
\end{figure}

Figure \ref{fig:results_gantt_all} illustrates a MA solution for a representative instance of each dataset, showing both the Gantt chart and the corresponding power requirement profile.
Expectedly, the algorithm consistently schedules jobs such that the resulting power requirements are aligned with moments where either there is sufficient on-site renewable generation, or the carbon intensity of the grid is at its lowest.
Depending on the available slack time in each instance, the algorithm may also introduce some planned idle time between operations.

\subsubsection{Algorithm performance}
\label{subsubsection:algorithm_performance}
\begin{figure}[h]
	\centering
		\begin{subfigure}{.49\columnwidth}
			\centering
			\includegraphics[clip, trim={0cm 0cm 0cm 0cm},width=\textwidth]{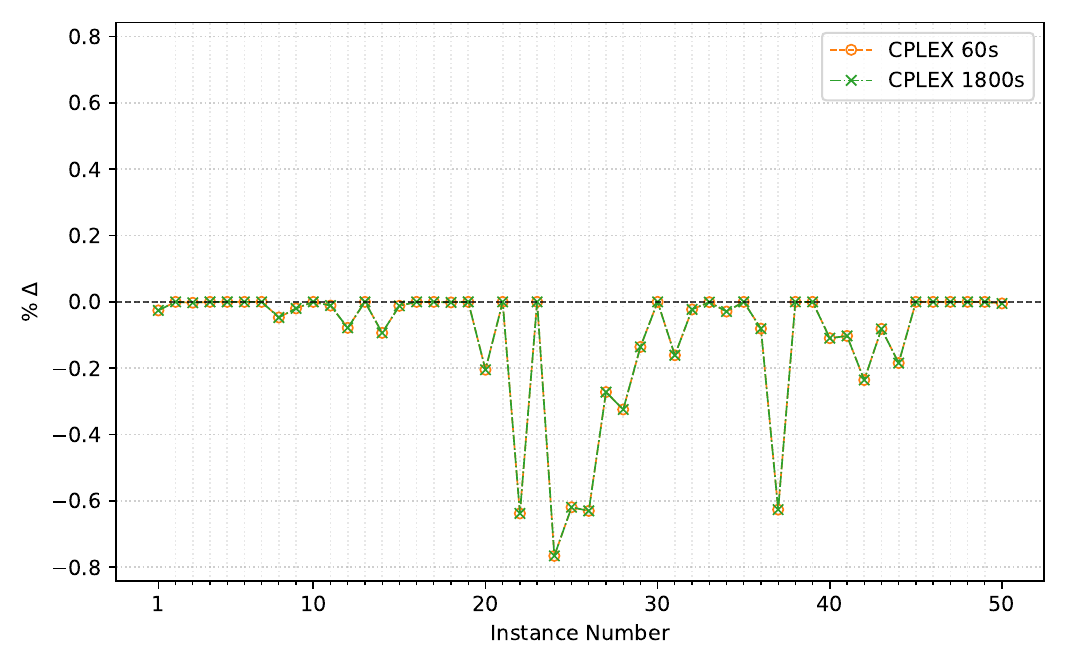}
			\caption{CAS-PFSP-M1T1 dataset}
			\label{fig:results_gap_per_instance_M1T1}
		\end{subfigure}
		\begin{subfigure}{.49\columnwidth}
			\centering
			\includegraphics[clip, trim={0cm 0cm 0cm 0cm},width=\textwidth]{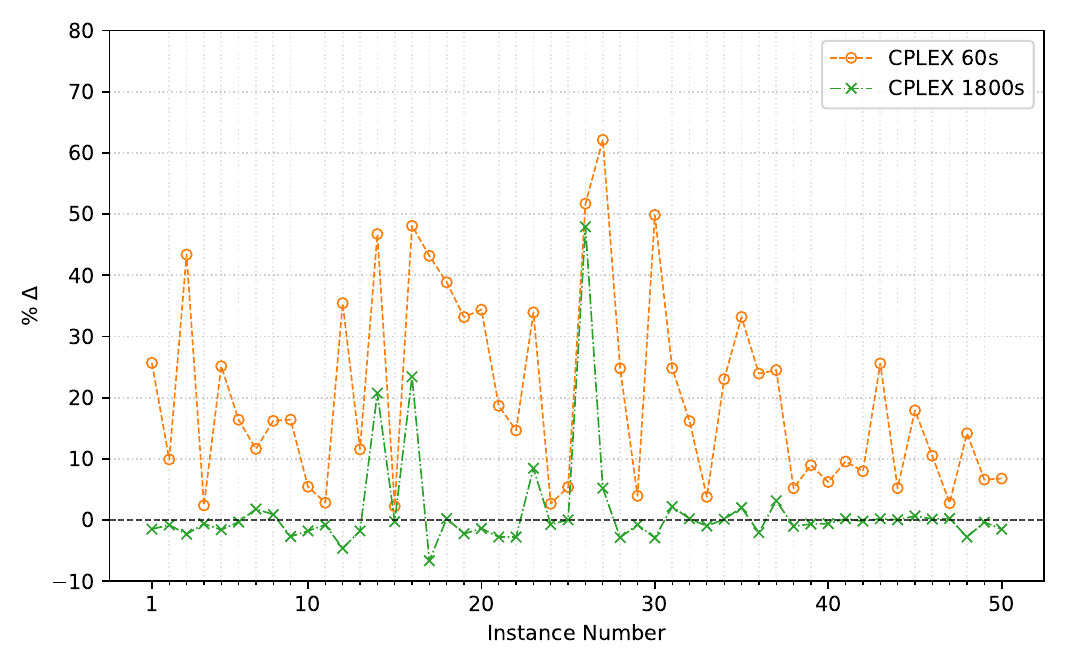}
			\caption{CAS-PFSP-M1T3 dataset}
			\label{fig:results_gap_per_instance_M1T3}
		\end{subfigure}
		\begin{subfigure}{.49\columnwidth}
			\centering
			\includegraphics[clip, trim={0cm 0cm 0cm 0cm},width=\textwidth]{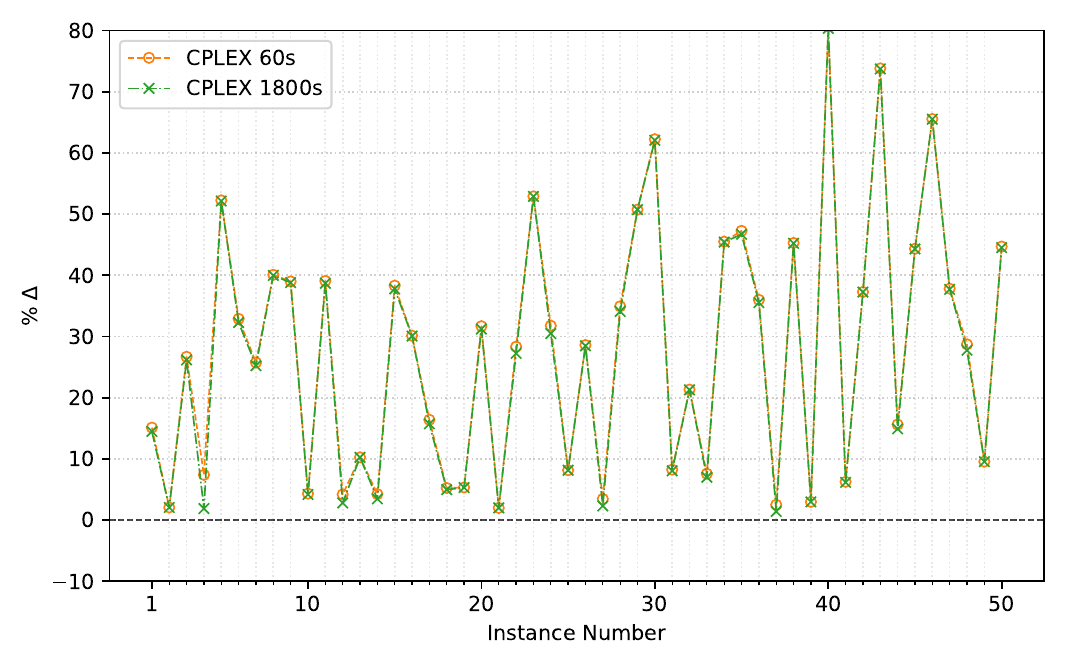}
			\caption{CAS-PFSP-M3T1 dataset}
			\label{fig:results_gap_per_instance_M3T1}
		\end{subfigure} 
		\begin{subfigure}{.49\columnwidth}
			\centering
			\includegraphics[clip, trim={0cm 0cm 0cm 0cm},width=\textwidth]{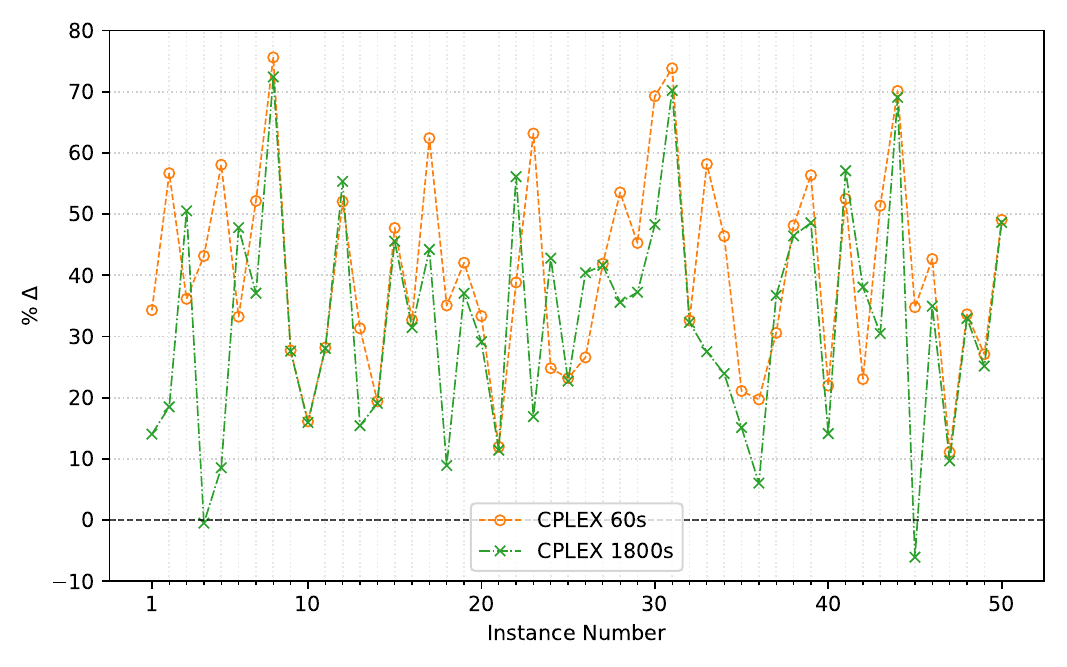}
			\caption{CAS-PFSP-M3T3 dataset}
			\label{fig:results_gap_per_instance_M3T3}
		\end{subfigure}
	\caption{Overview of percentage improvement per instance, per dataset}
	\label{fig:results_gap_per_instance_all}
\end{figure}

\begin{figure}[h]
	\centering
	\begin{subfigure}{.49\columnwidth}
		\includegraphics[width=\textwidth]{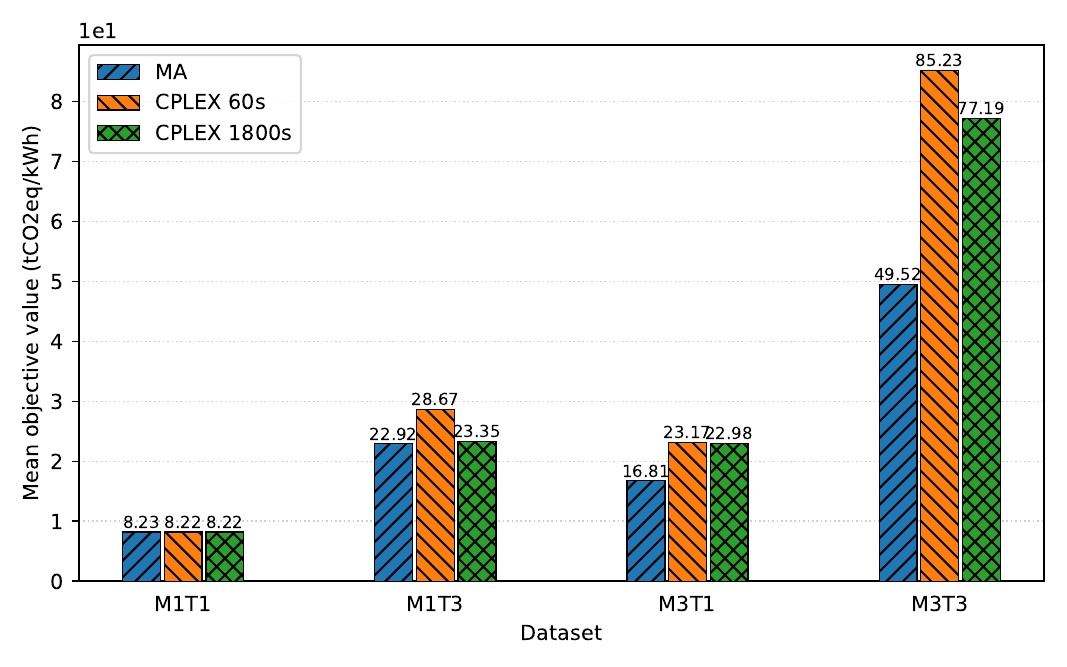}
		\caption{Objective value}
		\label{fig:results_bar_chart_obj}
	\end{subfigure}
	\begin{subfigure}{.49\columnwidth}
		\includegraphics[width=\textwidth]{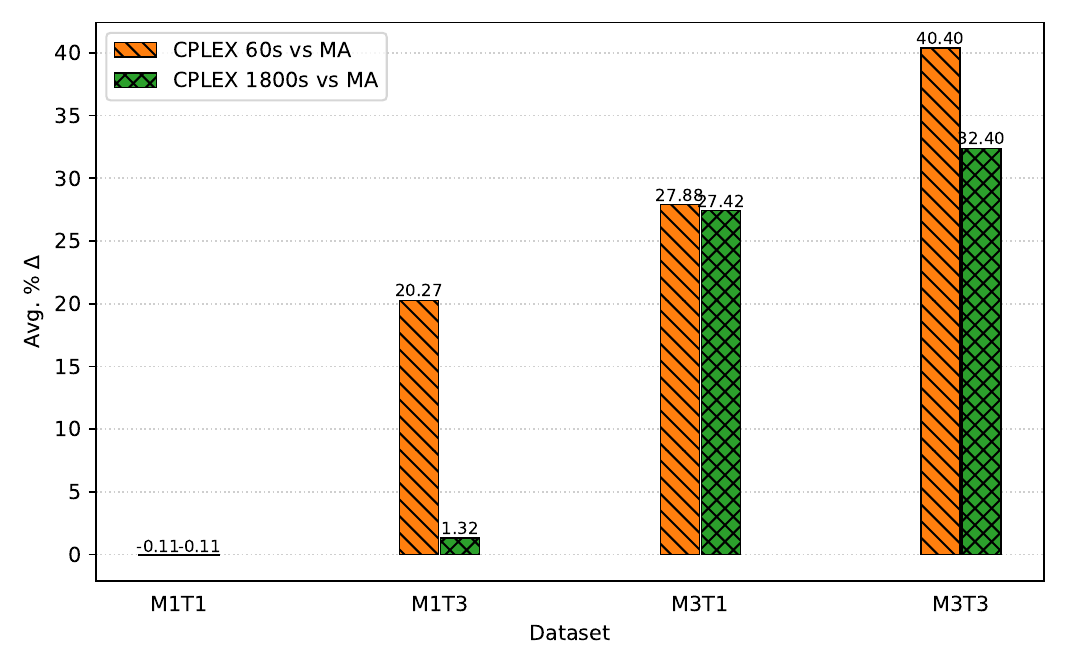}
		\caption{Percentage gap}
		\label{fig:results_bar_chart_gap}
	\end{subfigure}
	\caption{Average performance across each dataset per solution method}
	\label{fig:results_bar_chart_all}
\end{figure}

An overview of the MA's performance is reported in Figures \ref{fig:results_gap_per_instance_all} and \ref{fig:results_bar_chart_gap}, showing the percentage gaps across the four datasets both per instance and on average, respectively.

With both time limits of $1$ and $30$ minutes, the MA-CAS-PFSP algorithm consistently outperforms CPLEX on all datasets except CAS-PFSP-M1T1 (Figure \ref{fig:results_gap_per_instance_M1T1}).
This result was expected, as CPLEX found optimal solutions for all instances in this dataset, making it theoretically impossible for the MA to achieve better results.
Nevertheless, the MA obtains an average gap of $-0.11\, \%$ with an average computation time of $1.11$ seconds, compared to $4.45$ seconds for CPLEX.
These results confirm the effectiveness of the MA on smaller instances, which serves as a crucial validation step before applying it to larger instances where optimal solutions cannot be obtained with CPLEX.

The average percentage gap increases with the complexity of the instance dataset, as illustrated in Figure \ref{fig:results_bar_chart_gap}.
As expected, granting CPLEX more computational time improves its results across all datasets, compared to its own performance under the 1-minute time limit.
However, this improvement is more substantial on the longer time horizon datasets (Figures \ref{fig:results_gap_per_instance_M1T3} and \ref{fig:results_gap_per_instance_M3T3}) than on the shorter ones (Figures \ref{fig:results_gap_per_instance_M1T1} and \ref{fig:results_gap_per_instance_M3T1}).
This difference results from the combinatorial complexity of the instances in these datasets, which increases with the length of the time horizon.
In a PFSP scheduling environment, the sequence of the jobs remains fixed across all the machines.
Therefore, extending the time horizon (i.e., M1T3 compared to M1T1) has a greater impact on experimental problem complexity than increasing the number of machines (i.e., M3T1 compared to M1T1), which explains why additional computational time benefits CPLEX more on the longer datasets, although still not enough to outperform the MA.
Detailed per-instance results are available at \url{https://github.com/ugent-isye/CAS-PFSP}.

Figure \ref{fig:results_bar_chart_obj} presents the average objective value across all instances for each dataset.
Consistent with the previous discussion, performance gain in average objective value correlates with dataset complexity, highlighting the MA's strong scalability compared to CPLEX.
Average computational times are reported in Table \ref{tab:computational_times}.

\begin{table}[h]
	\fontsize{9pt}{14pt}\selectfont
	\centering
    \caption{Average computational times per dataset and solution method}
	\begin{tabular}{c|c|c|c}
		Dataset &  \multicolumn{3}{c}{Average time (s)} \\
		\cline{2-4}
		& MA & CPLEX 60s & CPLEX 1800s \\
		\hline
		CAS-MA-PFSP-M1T1 & 1.11 & 4.45 & 4.45 \\
		CAS-MA-PFSP-M1T3 & 2.87 & 60 & 1800 \\
		CAS-MA-PFSP-M3T1 & 7.59 & 60 & 1800 \\
		CAS-MA-PFSP-M3T3 & 24.49 & 60 & 1800 \\
	\end{tabular}
	\label{tab:computational_times}
\end{table}

\subsubsection{GHG emissions reduction}
\label{subsubsection:GHG_emissions_reduction}
To quantify the reduction in GHG emissions achieved with the proposed method, we compare the obtained solutions against two alternative optimization objectives: minimal makespan and minimal energy cost.

When optimizing for minimal makespan, the goal is to reduce the total completion time of the schedule while still respecting all precedence constraints.
In our case, the makespan is calculated as the completion time of the last operation executed on the final machine in the line.
This objective is among the most commonly used in scheduling literature and serves as a practical baseline for quantifying the emission reductions enabled by carbon-aware scheduling.

To provide additional insights, we also generate schedules optimized for minimal energy cost.
This is done by incorporating historical day-ahead electricity prices from the Belgian Spot BELPEX market.
These price data, available at \cite{ENTSO-E_electricity_prices}, represent hourly electricity tariffs and are matched to the energy-related data of each instance using their corresponding date.

For each of the three optimization targets (i.e., carbon minimization, cost minimization, and makespan minimization), we report the average objective values of the three objectives among the $10$ MA runs.
This allows us to assess the trade-offs between the objectives when optimizing for each one individually.

An overview of the results is presented in Table \ref{tab:overview_GHG_emissions}, where values are reported both in absolute terms and as relative percentage differences compared to the best value obtained for each objective (i.e., the row with the minimum value per column is taken as a reference and shown as $0.0 \, \%$).

In the absolute values section of the table, each row shows the average total GHG emissions, total energy cost, and makespan obtained when optimizing for the respective objective.
The relative values section reports the percentage increase over the best value achieved among the three approaches for each objective.
For example, in dataset CAS-PFSP-M1T1, on average, optimizing for carbon emissions leads to $0.0 \, \%$ additional emissions (by definition), $7.4 \, \%$ higher energy cost, and $5.5 \, \%$ longer makespan compared to the respective minima.

\begin{table}[h]
	\fontsize{9pt}{14pt}\selectfont
	\centering
    \caption{Overview of the results}
	\begin{tabular}{c|c||c|c|c||c|c|c}
        Dataset & Objective & \multicolumn{3}{c||}{Absolute values} & \multicolumn{3}{c}{Relative values} \\
        \cline{3-8}
        &  & Carbon (tCO\textsubscript{2}eq), & Cost (\euro{})& Makespan (min) & Carbon & Cost & Makespan\\
        \hline
                        & Carbon-min    & 8.23  & 5923  & 1413  & 0.0   & 7.4   & 5.5\\
        CAS-PFSP-M1T1   & Cost-min      & 8.49  & 5513  & 1432  & 3.2   & 0.0   & 6.9\\
                        & Makespan-min  & 10.51 & 7675  & 1339  & 27.7  & 39.2  & 0.0\\
        \hline
                        & Carbon-min    & 22.92 & 18049 & 4296  & 0.0   & 6.8   & 1.8\\
        CAS-PFSP-M1T3   & Cost-min      & 24.37 & 16894 & 4296  & 6.3   & 0.0   & 1.8\\
                        & Makespan-min  & 33.82 & 25203 & 4219  & 47.6  & 49.2  & 0.0\\
        \hline
                        & Carbon-min    & 16.81 & 11995  & 1251  & 0.0   & 11.7  & 44.5\\
        CAS-PFSP-M3T1   & Cost-min      & 17.34 & 10735  & 1228  & 3.2   & 0.0   & 41.8\\
                        & Makespan-min  & 21.66 & 15211  & 886   & 28.9  & 41.7  & 0.0 \\
        \hline
                        & Carbon-min    & 49.52 & 36353 & 4043  & 0.0   & 12.2  & 39.4\\
        CAS-PFSP-M3T3   & Cost-min      & 53.09 & 32404 & 4057  & 7.2   & 0.0   & 39.9\\
                        & Makespan-min  & 72.49 & 51159 & 2900  & 46.4  & 57.9  & 0.0\\
    \end{tabular}
    \label{tab:overview_GHG_emissions}
\end{table}

The results in Table \ref{tab:overview_GHG_emissions} reveal several key insights regarding the effects of different scheduling objectives on carbon emissions, energy cost, and makespan.
First, carbon-aware scheduling demonstrates a clear advantage in reducing GHG emissions across all datasets.
By leveraging time-dependent grid carbon intensity, on-site renewable electricity availability, and job-specific power requirements, it consistently achieves the lowest emissions levels.
Second, prioritizing makespan leads to a substantial increase in both emissions and energy costs.
Schedules optimized for minimal makespan result in $27.7 \, \%$ to $47.6 \, \%$ higher emissions, with energy costs rising even more sharply, from $39.2 \, \%$ to $57.9 \, \%$.
Finally, while minimizing energy costs does offer some environmental benefit compared to makespan reduction alone, it still results in considerably higher emissions than carbon-aware scheduling.
Across all datasets, cost-optimized schedules emit between $3.2 \, \%$ and $7.2 \, \%$ more GHGs than those optimized for carbon.

An important additional consideration is that optimizing for carbon emissions does not necessarily lead to a substantial increase in makespan.
This trade-off is particularly favorable in the single-machine datasets CAS-PFSP-M1T1 and CAS-PFSP-M1T3.
In the latter, the carbon-aware schedules achieve an average emission reduction of $47.6 \, \%$ with only a $1.8 \, \%$ increase in makespan.
By contrast, the multiple-machine datasets exhibit more pronounced makespan increases, peaking at $44.5 \, \%$.
This is largely due to the greater scheduling flexibility in these instances, which, as indicated in Table \ref{tab:overview_instance_datasets}, have considerably higher slack times.
Greater slack allows operations to be redistributed more freely across the planning horizon, which can result in longer makespans under carbon-aware scheduling.
Nonetheless, the single-machine results demonstrate that such flexibility is not a prerequisite for effective emission reductions: even with limited slack, substantial GHG reductions can be obtained with minimal impact on makespan.

\section{Conclusion}
\label{section:conclusion}
In this study, we have developed an efficient tool for aligning power consumption in manufacturing with the time-dependent carbon intensity of the grid and the availability of on-site renewable electricity, thereby reducing scope 2 GHG emissions.

The carbon-aware scheduling model was introduced in Section \ref{section:model_formulation}, where we first outlined how grid carbon intensity can be derived from from the generation mix and then formulated the scheduling problem as an MILP.
Given the NP-hardness of the problem, we proposed a dual random-key memetic algorithm in Section \ref{section:memetic_algorithm_framework} capable of efficiently finding high-quality solutions on real-world-sized instances by combining evolutionary computing with local search.
Finally, in Section \ref{section:computational_experiments}, we conducted computational experiments to evaluate the proposed approach by first discussing instance generation, then the performance evaluation metrics, and finally the experimental results.

We performed experiments by running the MA-CAS-PFSP algorithm and the CPLEX implementation of the MILP model on four instance datasets of increasing complexity.
The proposed algorithm outperformed CPLEX both in terms of solution quality and computation times on real-world sized instances.
To evaluate the benefits of carbon-aware scheduling, we compared it against two commonly used optimization objectives: makespan minimization and energy cost minimization.
The results show that carbon-aware scheduling consistently reduces GHG emissions across all datasets, achieving up to $47.6 \, \%$ lower scope 2 emissions compared to makespan-minimized schedules.
While energy costs minimization yields some environmental improvement over makespan minimization, it still results in significantly higher emissions, between $3.2 \, \%$ and $7.2 \, \%$ more, than carbon-aware scheduling.
Importantly, our analysis reveals that prioritizing GHG emissions in the objective function does not necessarily lead to a substantial increase in makespan.
In fact, in several cases, particularly in single-machine instances with limited scheduling flexibility, the increase in makespan remained marginal.

These findings highlight three key contributions of this study.
First, the study demonstrates how carbon-aware scheduling effectively reduces GHG emissions in manufacturing by integrating time-dependent grid's carbon intensity, on-site renewable generation, and job-specific power requirements within a PFSP scheduling model.
Second, the proposed algorithm efficiently solves large-scale instances, making it suitable for practical implementation.
Finally, the results show that scheduling for carbon emissions reduction does not necessarily require substantial sacrifices in makespan, particularly in scenarios with limited scheduling flexibility.

Despite these contributions, some limitations remain, offering directions for future work.
First, while the PFSP environment resembles many traditional assembly lines, it does not fully reflect the flexibility often required in modern manufacturing.
With the shift toward high-mix, low-volume production, manufacturers are adopting flexible assembly lines, where job routings can be dynamically adjusted to accommodate product customization.
Extending the model to a flexible-job-shop environment (FJSP) would enhance its applicability to contemporary manufacturing systems.
Second, the performance of metaheuristics is highly dependent on the values of their parameters, which require careful tuning to achieve good results across instances with similar characteristics.
However, energy-related data such as grid carbon intensity and on-site renewable generation can differ significantly from day to day, leading to substantial differences even between similar instances (e.g., same number of jobs, machines, and job power requirements).
These day-to-day differences may require frequent re-tuning of parameters, which is impractical for real-world deployment.
Moreover, the effectiveness of a parameter setting may change during the search.
Future research should explore self-adapting parameter control strategies, enabling the algorithm to dynamically adjust its parameters throughout the search process based on instance characteristics, thereby eliminating the need for manual (re-)tuning.
Third, while day-ahead forecasts for energy-related data are widely available, they are inherently uncertain.
Moreover, the accuracy of these forecasts diminishes as the time horizon extends.
As a result, a schedule that is initially optimal can become suboptimal as new, more accurate data becomes available throughout the planning horizon.
To mitigate this issue, future research should integrate real-time information into closed-loop rescheduling algorithms, allowing schedules to adapt dynamically and restore optimality as conditions change.

\section*{Acknowledgments}
We gratefully acknowledge the financial support provided by the Actemium chair on Sustainable Energy\footnote{\url{https://www.ugent.be/schenken/en/how-to-support/chairs/actemium}}.

\begingroup
    \setlength{\bibsep}{10pt}
    \setstretch{1}
	\bibliographystyle{apa} 
    \bibliography{refs_CAS}

\begin{thebibliography}{}

\bibitem[\protect\astroncite{Ahlqvist et~al.}{2022}]{Ahlqvist2022}
Ahlqvist, V., Holmberg, P., and Tangerås, T. (2022).
\newblock A survey comparing centralized and decentralized electricity markets.
\newblock {\em Energy Strategy Reviews}, 40:100812.
\newblock DOI:
  \href{http://doi.org/10.1016/j.esr.2022.100812}{10.1016/j.esr.2022.100812}.

\bibitem[\protect\astroncite{Ahmadian et~al.}{2021}]{Ahmadian2021}
Ahmadian, M.~M., Khatami, M., Salehipour, A., and Cheng, T. (2021).
\newblock Four decades of research on the open-shop scheduling problem to
  minimize the makespan.
\newblock {\em European Journal of Operational Research}, 295(2):399--426.
\newblock DOI:
  \href{http://doi.org/10.1016/j.ejor.2021.03.026}{10.1016/j.ejor.2021.03.026}.

\bibitem[\protect\astroncite{Akiba et~al.}{2019}]{Optuna2019}
Akiba, T., Sano, S., Yanase, T., Ohta, T., and Koyama, M. (2019).
\newblock Optuna: A next-generation hyperparameter optimization framework.
\newblock In {\em Proceedings of the 25th ACM SIGKDD International Conference
  on Knowledge Discovery \& Data Mining}, KDD '19, page 2623–2631, New York,
  NY, USA. Association for Computing Machinery.
\newblock DOI:
  \href{https://doi.org/10.1145/3292500.3330701}{10.1145/3292500.3330701}.

\bibitem[\protect\astroncite{Alidaee and Womer}{1999}]{Alidaee1999}
Alidaee, B. and Womer, N.~K. (1999).
\newblock Scheduling with time dependent processing times: Review and
  extensions.
\newblock {\em The Journal of the Operational Research Society},
  50(7):711--720.
\newblock DOI: \href{https://doi.org/10.2307/3010325}{10.2307/3010325}.

\bibitem[\protect\astroncite{Baker and Trietsch}{2018}]{Baker2018}
Baker, K.~R. and Trietsch, D. (2018).
\newblock {\em Principles of sequencing and scheduling}.
\newblock John Wiley \& Sons.
\newblock DOI:
  \href{http://doi.org/10.1002/9780470451793}{10.1002/9780470451793}.

\bibitem[\protect\astroncite{Batalla-Bejerano and
  Trujillo-Baute}{2016}]{Batalla-Bejerano2016}
Batalla-Bejerano, J. and Trujillo-Baute, E. (2016).
\newblock Impacts of intermittent renewable generation on electricity system
  costs.
\newblock {\em Energy Policy}, 94:411--420.
\newblock DOI:
  \href{http://doi.org/10.1016/j.enpol.2015.10.024}{10.1016/j.enpol.2015.10.024}.

\bibitem[\protect\astroncite{Baumgärtner et~al.}{2019}]{BAUMGARTNER2019755}
Baumgärtner, N., Delorme, R., Hennen, M., and Bardow, A. (2019).
\newblock Design of low-carbon utility systems: Exploiting time-dependent grid
  emissions for climate-friendly demand-side management.
\newblock {\em Applied Energy}, 247:755--765.
\newblock DOI:
  \href{http://doi.org/10.1016/j.apenergy.2019.04.029}{10.1016/j.apenergy.2019.04.029}.

\bibitem[\protect\astroncite{Bean}{1994}]{Bean1994}
Bean, J.~C. (1994).
\newblock Genetic algorithms and random keys for sequencing and optimization.
\newblock {\em ORSA journal on computing}, 6(2):154--160.
\newblock DOI:
  \href{https://doi.org/10.1287/ijoc.6.2.154}{10.1287/ijoc.6.2.154}.

\bibitem[\protect\astroncite{Bergstra et~al.}{2011}]{Bergstra2011}
Bergstra, J., Bardenet, R., Bengio, Y., and K\'{e}gl, B. (2011).
\newblock Algorithms for hyper-parameter optimization.
\newblock In {\em Advances in Neural Information Processing Systems},
  volume~24, pages 2546--2554.
\newblock
  \url{https://proceedings.neurips.cc/paper_files/paper/2011/file/86e8f7ab32cfd12577bc2619bc635690-Paper.pdf}.

\bibitem[\protect\astroncite{Che et~al.}{2016}]{Che2016}
Che, A., Zeng, Y., and Lyu, K. (2016).
\newblock An efficient greedy insertion heuristic for energy-conscious single
  machine scheduling problem under time-of-use electricity tariffs.
\newblock {\em Journal of Cleaner Production}, 129:565--577.
\newblock DOI:
  \href{https://doi.org/10.1016/j.jclepro.2016.03.150}{10.1016/j.jclepro.2016.03.150}.

\bibitem[\protect\astroncite{Cheng et~al.}{2017}]{Cheng2017}
Cheng, J., Chu, F., Liu, M., Wu, P., and Xia, W. (2017).
\newblock Bi-criteria single-machine batch scheduling with machine on/off
  switching under time-of-use tariffs.
\newblock {\em Computers \& Industrial Engineering}, 112:721--734.
\newblock DOI:
  \href{https://doi.org/10.1016/j.cie.2017.04.026}{10.1016/j.cie.2017.04.026}.

\bibitem[\protect\astroncite{Cheng et~al.}{1999}]{Cheng1999}
Cheng, R., Gen, M., and Tsujimura, Y. (1999).
\newblock A tutorial survey of job-shop scheduling problems using genetic
  algorithms, part ii: hybrid genetic search strategies.
\newblock {\em Computers \& Industrial Engineering}, 36(2):343--364.
\newblock DOI:
  \href{http://doi.org/10.1016/S0360-8352(99)00136-9}{10.1016/S0360-8352(99)00136-9}.

\bibitem[\protect\astroncite{Cheng et~al.}{2004}]{Cheng2004}
Cheng, T., Ding, Q., and Lin, B. (2004).
\newblock A concise survey of scheduling with time-dependent processing times.
\newblock {\em European Journal of Operational Research}, 152(1):1--13.
\newblock DOI:
  \href{https://doi.org/10.1016/S0377-2217(02)00909-8}{10.1016/S0377-2217(02)00909-8}.

\bibitem[\protect\astroncite{Cui and Lu}{2021}]{Cui2021}
Cui, W. and Lu, B. (2021).
\newblock Energy-aware operations management for flow shops under tou
  electricity tariff.
\newblock {\em Computers \& Industrial Engineering}, 151:106942.
\newblock DOI:
  \href{https://doi.org/10.1016/j.cie.2020.106942}{10.1016/j.cie.2020.106942}.

\bibitem[\protect\astroncite{da~Jiang and Wang}{2020}]{Jiang2020}
da~Jiang, E. and Wang, L. (2020).
\newblock Multi-objective optimization based on decomposition for flexible job
  shop scheduling under time-of-use electricity prices.
\newblock {\em Knowledge-Based Systems}, 204:106177.
\newblock DOI:
  \href{https://doi.org/10.1016/j.knosys.2020.106177}{10.1016/j.knosys.2020.106177}.

\bibitem[\protect\astroncite{Ding et~al.}{2016a}]{Ding2016Carbon}
Ding, J.-Y., Song, S., and Wu, C. (2016a).
\newblock Carbon-efficient scheduling of flow shops by multi-objective
  optimization.
\newblock {\em European Journal of Operational Research}, 248(3):758--771.
\newblock DOI:
  \href{https://doi.org/10.1016/j.ejor.2015.05.019}{10.1016/j.ejor.2015.05.019}.

\bibitem[\protect\astroncite{Ding et~al.}{2016b}]{Ding2016Parallel}
Ding, J.-Y., Song, S., Zhang, R., Chiong, R., and Wu, C. (2016b).
\newblock Parallel machine scheduling under time-of-use electricity prices: New
  models and optimization approaches.
\newblock {\em IEEE Transactions on Automation Science and Engineering},
  13(2):1138--1154.
\newblock DOI:
  \href{https://doi.org/10.1109/TASE.2015.2495328}{10.1109/TASE.2015.2495328}.

\bibitem[\protect\astroncite{Dong and Ye}{2022}]{Dong2022}
Dong, J. and Ye, C. (2022).
\newblock Green scheduling of distributed two-stage reentrant hybrid flow shop
  considering distributed energy resources and energy storage system.
\newblock {\em Computers \& Industrial Engineering}, 169:108146.
\newblock DOI:
  \href{http://doi.org/10.1016/j.cie.2022.108146}{10.1016/j.cie.2022.108146}.

\bibitem[\protect\astroncite{Eiben and Smith}{2015}]{eiben2015introduction}
Eiben, A.~E. and Smith, J.~E. (2015).
\newblock {\em Introduction to evolutionary computing}.
\newblock Springer.
\newblock DOI:
  \href{http://doi.org/10.1007/978-3-662-44874-8}{10.1007/978-3-662-44874-8}.

\bibitem[\protect\astroncite{Elia}{2025a}]{Elia_DataSet_GridMix}
Elia (2025a).
\newblock {\em Day-ahead schedule of power generation - aggregated by fuel type
  - data up to 22/05/2024}.
\newblock Open Data Platform.
\newblock \url{https://opendata.elia.be/explore/dataset/ods034/information/}
  Accessed: 2025-01-08.

\bibitem[\protect\astroncite{Elia}{2025b}]{Elia_DataSet_Onsite}
Elia (2025b).
\newblock {\em Photovoltaic power production estimation and forecast on Belgian
  grid (Historical)}.
\newblock Open Data Platform.
\newblock \url{https://opendata.elia.be/explore/dataset/ods032/information/}
  Accessed: 2025-01-08.

\bibitem[\protect\astroncite{{ENTSO-E}}{2025}]{ENTSO-E_electricity_prices}
{ENTSO-E} (2025).
\newblock {\em Central collection and publication of electricity generation,
  transportation and consumption data and information for the pan-European
  market}.
\newblock Transparency Platform.
\newblock \url{https://newtransparency.entsoe.eu/} Accessed: 2025-06-11.

\bibitem[\protect\astroncite{{European Commission and Directorate-General for
  Climate Action}}{2019}]{net-zero_plan_2019}
{European Commission and Directorate-General for Climate Action} (2019).
\newblock {\em Going climate-neutral by 2050 – A strategic long-term vision
  for a prosperous, modern, competitive and climate-neutral EU economy}.
\newblock Publications Office of the European Union.
\newblock DOI: \href{http://doi.org/10.2834/02074}{10.2834/02074}.

\bibitem[\protect\astroncite{{European Commission:
  Eurostat}}{tion}]{Eurostat_energy_report_2023}
{European Commission: Eurostat} (2023 interactive publication).
\newblock {\em Shredding light on energy in the {EU}}.
\newblock Publications Office of the European Union.
\newblock DOI: \href{http://doi.org/10.2785/405482}{10.2785/405482}.

\bibitem[\protect\astroncite{{European Environment
  Agency}}{2023}]{EEA_report_2023}
{European Environment Agency} (2023).
\newblock {\em Trends and projections in Europe 2023}.
\newblock Publications Office of the European Union.
\newblock DOI: \href{http://doi.org/10.2800/595102}{10.2800/595102}.

\bibitem[\protect\astroncite{{European Union: European
  Commission}}{inal}]{eu_greendeal2023}
{European Union: European Commission} (1 February 2023, COM(2023) 62 final).
\newblock {\em Communication from the Commission to the European Parliament,
  the European Council, the Council, the European Economic and Social Committee
  and the Committee of the Regions: A Green Deal Industrial Plan for the
  Net-Zero Age}.
\newblock The European Commission.
\newblock
  \url{https://eur-lex.europa.eu/legal-content/EN/TXT/?uri=CELEX%3A52023DC0062&qid=1678874483913}.

\bibitem[\protect\astroncite{Fallahi et~al.}{2023}]{Fallahi2023}
Fallahi, A., Shahidi-Zadeh, B., and Niaki, S. T.~A. (2023).
\newblock Unrelated parallel batch processing machine scheduling for production
  systems under carbon reduction policies: Nsga-ii and mogwo metaheuristics.
\newblock {\em Soft Computing}, 27(22):17063--17091.
\newblock DOI:
  \href{https://doi.org/10.1007/s00500-023-08754-0}{10.1007/s00500-023-08754-0}.

\bibitem[\protect\astroncite{Fazli~Khalaf and Wang}{2018}]{FazliKhalaf2018}
Fazli~Khalaf, A. and Wang, Y. (2018).
\newblock Energy-cost-aware flow shop scheduling considering intermittent
  renewables, energy storage, and real-time electricity pricing.
\newblock {\em International Journal of Energy Research}, 42(12):3928--3942.
\newblock DOI: \href{https://doi.org/10.1002/er.4130}{10.1002/er.4130}.

\bibitem[\protect\astroncite{Fernandez-Viagas
  et~al.}{2017}]{Fernandez-Viagas2017}
Fernandez-Viagas, V., Ruiz, R., and Framinan, J.~M. (2017).
\newblock A new vision of approximate methods for the permutation flowshop to
  minimise makespan: State-of-the-art and computational evaluation.
\newblock {\em European Journal of Operational Research}, 257(3):707--721.
\newblock DOI:
  \href{https://doi.org/10.1016/j.ejor.2016.09.055}{10.1016/j.ejor.2016.09.055}.

\bibitem[\protect\astroncite{Gabrek and Seifermann}{2025}]{Gabrek2025}
Gabrek, N. and Seifermann, S. (2025).
\newblock How the correlation between electricity prices and emission intensity
  affects the economic and ecological potential of industrial demand-side
  flexibility measures.
\newblock {\em Journal of Cleaner Production}, page 145863.
\newblock DOI:
  \href{https://doi.org/10.1016/j.jclepro.2025.145863}{10.1016/j.jclepro.2025.145863}.

\bibitem[\protect\astroncite{Garey et~al.}{1976}]{Garey1976}
Garey, M.~R., Johnson, D.~S., and Sethi, R. (1976).
\newblock The complexity of flowshop and jobshop scheduling.
\newblock {\em Mathematics of Operations Research}, 1(2):117--129.
\newblock DOI:
  \href{http://doi.org/10.1287/moor.1.2.117}{10.1287/moor.1.2.117}.

\bibitem[\protect\astroncite{{Ghorbani Saber} and
  Ranjbar}{2022}]{GhorbaniSaber2022}
{Ghorbani Saber}, R. and Ranjbar, M. (2022).
\newblock Minimizing the total tardiness and the total carbon emissions in the
  permutation flow shop scheduling problem.
\newblock {\em Computers \& Operations Research}, 138:105604.
\newblock DOI:
  \href{https://doi.org/10.1016/j.cor.2021.105604}{10.1016/j.cor.2021.105604}.

\bibitem[\protect\astroncite{Ghorbanzadeh and Ranjbar}{2023}]{Ghorbanzadeh2023}
Ghorbanzadeh, M. and Ranjbar, M. (2023).
\newblock Energy-aware production scheduling in the flow shop environment under
  sequence-dependent setup times, group scheduling and renewable energy
  constraints.
\newblock {\em European Journal of Operational Research}, 307(2):519--537.
\newblock DOI:
  \href{http://doi.org/10.1016/j.ejor.2022.09.034}{10.1016/j.ejor.2022.09.034}.

\bibitem[\protect\astroncite{{Greenhouse Gas
  Protocol}}{2004}]{GHGCorporateStandard}
{Greenhouse Gas Protocol} (2004).
\newblock {\em A Corporate Accounting and Reporting Standard}.
\newblock World Resources Institute and World Business Council for Sustainable
  Development.
\newblock
  \url{https://ghgprotocol.org/sites/default/files/standards/ghg-protocol-revised.pdf}.

\bibitem[\protect\astroncite{{Greenhouse Gas
  Protocol}}{2015}]{GHGscope2guidance}
{Greenhouse Gas Protocol} (2015).
\newblock {\em {GHG} {P}rotocol {S}cope 2 {G}uidance}.
\newblock WRI: World Resources Institute.
\newblock
  \url{https://ghgprotocol.org/sites/default/files/2023-03/Scope%202%20Guidance.pdf}.

\bibitem[\protect\astroncite{Ho et~al.}{2022}]{Ho2022}
Ho, M.~H., Hnaien, F., and Dugardin, F. (2022).
\newblock Exact method to optimize the total electricity cost in two-machine
  permutation flow shop scheduling problem under time-of-use tariff.
\newblock {\em Computers \& Operations Research}, 144:105788.
\newblock DOI:
  \href{https://doi.org/10.1016/j.cor.2022.105788}{10.1016/j.cor.2022.105788}.

\bibitem[\protect\astroncite{Holland and Mansur}{2008}]{Holland2008}
Holland, S.~P. and Mansur, E.~T. (2008).
\newblock {Is Real-Time Pricing Green? The Environmental Impacts of Electricity
  Demand Variance}.
\newblock {\em The Review of Economics and Statistics}, 90(3):550--561.
\newblock DOI:
  \href{http://doi.org/10.1162/rest.90.3.550}{10.1162/rest.90.3.550}.

\bibitem[\protect\astroncite{{IEA}}{2022}]{IEA2022}
{IEA} (2022).
\newblock {\em Electricity Market Report - January 2022}.
\newblock International Energy Agency, Paris.
\newblock
  \url{https://www.iea.org/reports/electricity-market-report-january-2022}.

\bibitem[\protect\astroncite{{Intergovernmental Panel on Climate
  Change}}{2015}]{IPCC_ar5}
{Intergovernmental Panel on Climate Change} (2015).
\newblock Technology-specific cost and performance parameters.
\newblock In {\em Climate Change 2014: Mitigation of Climate Change: Working
  Group III Contribution to the IPCC Fifth Assessment Report}, page
  1329–1356. Cambridge University Press.
\newblock DOI:
  \href{http://doi.org/10.1017/CBO9781107415416}{10.1017/CBO9781107415416}.

\bibitem[\protect\astroncite{{Intergovernmental Panel on Climate
  Change}}{2023}]{IPCC_AR6}
{Intergovernmental Panel on Climate Change} (2023).
\newblock {\em Climate Change 2023: Synthesis Report. Contribution of Working
  Groups I, II and III to the Sixth Assessment Report of the Intergovernmental
  Panel on Climate Change}.
\newblock IPCC, Geneva, Switzerland, pp. 35-115.
\newblock DOI:
  \href{http://doi.org/10.59327/IPCC/AR6-9789291691647}{10.59327/IPCC/AR6-9789291691647}.

\bibitem[\protect\astroncite{Jaehn and Sedding}{2016}]{Jaehn2016}
Jaehn, F. and Sedding, H.~A. (2016).
\newblock Scheduling with time-dependent discrepancy times.
\newblock {\em Journal of Scheduling}, 19(6):737--757.
\newblock DOI:
  \href{https://doi.org/10.1007/s10951-016-0472-2}{10.1007/s10951-016-0472-2}.

\bibitem[\protect\astroncite{Karimi and Kwon}{2021}]{Karimi2021}
Karimi, S. and Kwon, S. (2021).
\newblock Comparative analysis of the impact of energy-aware scheduling,
  renewable energy generation, and battery energy storage on production
  scheduling.
\newblock {\em International Journal of Energy Research}, 45(13):18981--18998.
\newblock DOI: \href{https://doi.org/10.1002/er.6999}{10.1002/er.6999}.

\bibitem[\protect\astroncite{Kelley et~al.}{2018}]{Kelley2018}
Kelley, M.~T., Baldick, R., and Baldea, M. (2018).
\newblock Demand response operation of electricity-intensive chemical processes
  for reduced greenhouse gas emissions: application to an air separation unit.
\newblock {\em ACS Sustainable Chemistry \& Engineering}, 7(2):1909--1922.
\newblock DOI:
  \href{https://doi.org/10.1021/acssuschemeng.8b03927}{10.1021/acssuschemeng.8b03927}.

\bibitem[\protect\astroncite{Kopsakangas-Savolainen
  et~al.}{2017}]{KOPSAKANGASSAVOLAINEN2017384}
Kopsakangas-Savolainen, M., Mattinen, M.~K., Manninen, K., and Nissinen, A.
  (2017).
\newblock Hourly-based greenhouse gas emissions of electricity – cases
  demonstrating possibilities for households and companies to decrease their
  emissions.
\newblock {\em Journal of Cleaner Production}, 153:384--396.
\newblock DOI:
  \href{http://doi.org/j.jclepro.2015.11.027}{j.jclepro.2015.11.027}.

\bibitem[\protect\astroncite{Liu}{2016}]{Liu2016}
Liu, C.-H. (2016).
\newblock Mathematical programming formulations for single-machine scheduling
  problems while considering renewable energy uncertainty.
\newblock {\em International Journal of Production Research}, 54(4):1122--1133.
\newblock DOI:
  \href{https://doi.org/10.1080/00207543.2015.1048380}{10.1080/00207543.2015.1048380}.

\bibitem[\protect\astroncite{Liu et~al.}{2017}]{Liu2017}
Liu, Q., Zhan, M., Chekem, F.~O., Shao, X., Ying, B., and Sutherland, J.~W.
  (2017).
\newblock A hybrid fruit fly algorithm for solving flexible job-shop scheduling
  to reduce manufacturing carbon footprint.
\newblock {\em Journal of Cleaner Production}, 168:668--678.
\newblock DOI:
  \href{https://doi.org/10.1016/j.jclepro.2017.09.037}{10.1016/j.jclepro.2017.09.037}.

\bibitem[\protect\astroncite{Londe et~al.}{2025}]{Londe2025}
Londe, M.~A., Pessoa, L.~S., Andrade, C.~E., and Resende, M.~G. (2025).
\newblock Biased random-key genetic algorithms: A review.
\newblock {\em European Journal of Operational Research}, 321(1):1--22.
\newblock DOI:
  \href{https://doi.org/10.1016/j.ejor.2024.03.030}{10.1016/j.ejor.2024.03.030}.

\bibitem[\protect\astroncite{Luo et~al.}{2013}]{Luo2013}
Luo, H., Du, B., Huang, G.~Q., Chen, H., and Li, X. (2013).
\newblock Hybrid flow shop scheduling considering machine electricity
  consumption cost.
\newblock {\em International Journal of Production Economics}, 146(2):423--439.
\newblock DOI:
  \href{https://doi.org/10.1016/j.ijpe.2013.01.028}{10.1016/j.ijpe.2013.01.028}.

\bibitem[\protect\astroncite{Miller et~al.}{2022}]{Miller2022}
Miller, G.~J., Novan, K., and Jenn, A. (2022).
\newblock Hourly accounting of carbon emissions from electricity consumption.
\newblock {\em Environmental Research Letters}, 17(4).
\newblock DOI:
  \href{http://doi.org/10.1088/1748-9326/ac6147}{10.1088/1748-9326/ac6147}.

\bibitem[\protect\astroncite{Minh Hung~Ho and Dugardin}{2021}]{Ho2021}
Minh Hung~Ho, F.~H. and Dugardin, F. (2021).
\newblock Electricity cost minimisation for optimal makespan solution in flow
  shop scheduling under time-of-use tariffs.
\newblock {\em International Journal of Production Research}, 59(4):1041--1067.
\newblock DOI:
  \href{https://doi.org/10.1080/00207543.2020.1715504}{10.1080/00207543.2020.1715504}.

\bibitem[\protect\astroncite{Moon and Park}{2014}]{Moon2014}
Moon, J.-Y. and Park, J. (2014).
\newblock Smart production scheduling with time-dependent and machine-dependent
  electricity cost by considering distributed energy resources and energy
  storage.
\newblock {\em International Journal of Production Research},
  52(13):3922--3939.
\newblock DOI:
  \href{https://doi.org/10.1080/00207543.2013.860251}{10.1080/00207543.2013.860251}.

\bibitem[\protect\astroncite{Moon et~al.}{2013}]{Moon2013}
Moon, J.-Y., Shin, K., and Park, J. (2013).
\newblock Optimization of production scheduling with time-dependent and
  machine-dependent electricity cost for industrial energy efficiency.
\newblock {\em The International Journal of Advanced Manufacturing Technology},
  68:523--535.
\newblock DOI:
  \href{https://doi.org/10.1007/s00170-013-4749-8}{10.1007/s00170-013-4749-8}.

\bibitem[\protect\astroncite{Ostermeier and Deuse}{2024}]{Ostermeier2024}
Ostermeier, F.~F. and Deuse, J. (2024).
\newblock A review and classification of scheduling objectives in unpaced flow
  shops for discrete manufacturing.
\newblock {\em Journal of Scheduling}, 27(1):29--49.
\newblock DOI:
  \href{http://doi.org/10.1007/s10951-023-00795-5}{10.1007/s10951-023-00795-5}.

\bibitem[\protect\astroncite{Oukil et~al.}{2022}]{Oukil2022}
Oukil, A., El-Bouri, A., and Emrouznejad, A. (2022).
\newblock Energy-aware job scheduling in a multi-objective production
  environment – an integrated dea-owa model.
\newblock {\em Computers \& Industrial Engineering}, 168:108065.
\newblock DOI:
  \href{https://doi.org/10.1016/j.cie.2022.108065}{10.1016/j.cie.2022.108065}.

\bibitem[\protect\astroncite{Pinedo}{2022}]{PinedoScheduling}
Pinedo, M.~L. (2022).
\newblock {\em Scheduling: Theory, Algorithms, and Systems}.
\newblock Springer Cham.
\newblock DOI:
  \href{http://doi.org/10.1007/978-3-031-05921-6}{10.1007/978-3-031-05921-6}.

\bibitem[\protect\astroncite{Rubaiee and Yildirim}{2019}]{Rubaiee2019}
Rubaiee, S. and Yildirim, M.~B. (2019).
\newblock An energy-aware multiobjective ant colony algorithm to minimize total
  completion time and energy cost on a single-machine preemptive scheduling.
\newblock {\em Computers \& Industrial Engineering}, 127:240--252.
\newblock DOI:
  \href{https://doi.org/10.1016/j.cie.2018.12.020}{10.1016/j.cie.2018.12.020}.

\bibitem[\protect\astroncite{Shen et~al.}{2021}]{Shen2021}
Shen, K., {De Pessemier}, T., Martens, L., and Joseph, W. (2021).
\newblock A parallel genetic algorithm for multi-objective flexible flowshop
  scheduling in pasta manufacturing.
\newblock {\em Computers \& Industrial Engineering}, 161:107659.
\newblock DOI:
  \href{https://doi.org/10.1016/j.cie.2021.107659}{10.1016/j.cie.2021.107659}.

\bibitem[\protect\astroncite{Shen et~al.}{2023}]{Shen2023}
Shen, L., Dauzère-Pérès, S., and Maecker, S. (2023).
\newblock Energy cost efficient scheduling in flexible job-shop manufacturing
  systems.
\newblock {\em European Journal of Operational Research}, 310(3):992--1016.
\newblock DOI:
  \href{http://doi.org/10.1016/j.ejor.2023.03.041}{10.1016/j.ejor.2023.03.041}.

\bibitem[\protect\astroncite{Soares and Carvalho}{2020}]{Soares2020}
Soares, L. C.~R. and Carvalho, M. A.~M. (2020).
\newblock Biased random-key genetic algorithm for scheduling identical parallel
  machines with tooling constraints.
\newblock {\em European Journal of Operational Research}, 285(3):955--964.
\newblock DOI:
  \href{https://doi.org/10.1016/j.ejor.2020.02.047}{10.1016/j.ejor.2020.02.047}.

\bibitem[\protect\astroncite{Stecco et~al.}{2008}]{Stecco2008}
Stecco, G., Cordeau, J.-F., and Moretti, E. (2008).
\newblock A branch-and-cut algorithm for a production scheduling problem with
  sequence-dependent and time-dependent setup times.
\newblock {\em Computers \& Operations Research}, 35(8):2635--2655.
\newblock DOI:
  \href{https://doi.org/10.1016/j.cor.2006.12.021}{10.1016/j.cor.2006.12.021}.

\bibitem[\protect\astroncite{Tian and Zheng}{2024}]{Tian2024}
Tian, Z. and Zheng, L. (2024).
\newblock Single machine parallel-batch scheduling under time-of-use
  electricity prices: New formulations and optimisation approaches.
\newblock {\em European Journal of Operational Research}, 312(2):512--524.
\newblock DOI:
  \href{https://doi.org/10.1016/j.ejor.2023.07.012}{10.1016/j.ejor.2023.07.012}.

\bibitem[\protect\astroncite{Trevino-Martinez
  et~al.}{2022a}]{Trevinomartinez2022Footprintoptimization}
Trevino-Martinez, S., Sawhney, R., and Shylo, O. (2022a).
\newblock Energy-carbon footprint optimization in sequence-dependent production
  scheduling.
\newblock {\em Applied Energy}, 315:118949.
\newblock DOI:
  \href{https://doi.org/10.1016/j.apenergy.2022.118949}{10.1016/j.apenergy.2022.118949}.

\bibitem[\protect\astroncite{Trevino-Martinez
  et~al.}{2022b}]{Trevinomartinez2022Neutralityoptimization}
Trevino-Martinez, S., Sawhney, R., and Sims, C. (2022b).
\newblock Energy-carbon neutrality optimization in production scheduling via
  solar net metering.
\newblock {\em Journal of Cleaner Production}, 380:134627.
\newblock DOI:
  \href{https://doi.org/10.1016/j.jclepro.2022.134627}{10.1016/j.jclepro.2022.134627}.

\bibitem[\protect\astroncite{Wang et~al.}{2018}]{Wang2018}
Wang, S., Zhu, Z., Fang, K., Chu, F., and Chu, C. (2018).
\newblock Scheduling on a two-machine permutation flow shop under time-of-use
  electricity tariffs.
\newblock {\em International Journal of Production Research}, 56(9):3173--3187.
\newblock DOI:
  \href{https://doi.org/10.1080/00207543.2017.1401236}{10.1080/00207543.2017.1401236}.

\bibitem[\protect\astroncite{Wolpert and Macready}{1997}]{Wolpert1997}
Wolpert, D. and Macready, W. (1997).
\newblock No free lunch theorems for optimization.
\newblock {\em {IEEE} Transactions on Evolutionary Computation}, 1(1):67--82.
\newblock DOI: \href{https://doi.org/10.1109/4235.585893}{10.1109/4235.585893}.

\bibitem[\protect\astroncite{Zhang et~al.}{2014}]{Zhang2014}
Zhang, H., Zhao, F., Fang, K., and Sutherland, J.~W. (2014).
\newblock Energy-conscious flow shop scheduling under time-of-use electricity
  tariffs.
\newblock {\em CIRP Annals}, 63(1):37--40.
\newblock DOI:
  \href{https://doi.org/10.1016/j.cirp.2014.03.011}{10.1016/j.cirp.2014.03.011}.

\end{thebibliography}
\endgroup

\end{document}